\documentclass[11pt,reqno]{amsart} % for final version

\newcommand\ifshort\iftrue % for article
\newcommand\iflong\iffalse

\usepackage{bbm, amsmath, amsfonts, amscd, latexsym, amsthm, amssymb, graphicx,hyperref,mathrsfs, comment}%, refcheck}
\usepackage{abstract}
\usepackage{mathtools}
\usepackage{caption}
\usepackage[dvipsnames]{xcolor}

\usepackage{tikz-cd}
\tikzset{
  symbol/.style={
    draw=none,
    every to/.append style={
      edge node={node [sloped, allow upside down, auto=false]{$#1$}}}
  }
}
\usepackage[shortalphabetic]{amsrefs}

\makeatletter
\newif\if@check@engine  \@check@enginetrue 
\makeatother
\usepackage[usenames,dvipsnames]{pstricks}

\newtheorem{theor}{\hspace{1cm}{\sc Theorem}}[section]
\newtheorem{utver}[theor]{\hspace{1cm}{\sc Proposition}}
\newtheorem{sledst}[theor]{\hspace{1cm}{\sc Corollary}}
\newtheorem{lemma}[theor]{\hspace{1cm}{\sc Lemma}}

\newtheorem*{utver*}{\hspace{1cm}{\sc Proposition}}
\theoremstyle{definition}

\newtheorem{defin}[theor]{\hspace{1cm}{\sc Definition}}
\newtheorem{exa}[theor]{\hspace{1cm}{\sc Example}}

\newtheorem{rem}[theor]{\hspace{1cm}{\sc Remark}}

\newcommand{\codim}{\mathop{\rm codim}\nolimits}

\newcommand{\sing}{\mathop{\rm sing}\nolimits}

\newcommand{\cork}{\mathop{\rm cork}\nolimits}

\newcommand{\conv}{\mathop{\rm conv}\nolimits}

\renewcommand{\Im}{\mathop{\rm Im}\nolimits}

\newcommand{\ord}{\mathop{\rm ord}\nolimits}

\renewcommand{\emph}[1]{{\it {\color{NavyBlue} #1}}}

\def\R{\mathbb R}
\def\N{\mathbb N}
\def\Z{\mathbb Z}

\def\C{\mathbb C}
\def\CC{({\mathbb C}^\star)}

\def\CP{\mathbb C\mathbb P}

\def\Pr{\mathbb P}

\def\x{\times}

\title{Singular loci of sparse resultants}
\author{Evgeny Statnik$^1$}
\address{$^1$\rm HSE University}\address{\rm
Moscow, Russia}\address{\tt evstatnik@gmail.com}
\address{\rm The author is partially supported by International Laboratory of Cluster Geometry NRU HSE, RF Government grant, ag. № 075-15-2021-608 dated 08.06.2021}

\usepackage[foot]{amsaddr}

\begin{document}
\maketitle

\begin{abstract}
We study the singularity locus of the sparse resultant of two univariate polynomials, and apply our results to estimate singularities of a coordinate projection of a generic spatial complete intersection curve.
\end{abstract}

\tableofcontents

\section{Introduction}

Given a pair of finite sets $B_1$ and $B_2\subset\Z$ with at least two elements in each, one can consider the space of pairs of sparse Laurent polynomials supported at $B_1$ and $B_2$:
$$\C^{B_1}\times\C^{B_2}=\left\{\left(f_1(x)=\sum_{k\in B_1}f^1_kx^k,\,f_2(x)=\sum_{k\in B_2}f^2_kx^k\right)\right\}.$$

The sparse resultant $R_B\subset\C^{B_1} \times \C^{B_2}$ is the closure of the set
$$H(1)=\{(f_1,f_2)\in\C^{B_1} \times \C^{B_2}\,|\,f_1(x)=f_2(x)=0\mbox{ has at least one non-zero solution in } \C^*\}.$$

\begin{rem}
We need the closure here, because the roots can tend to $0$ and infinity. % It takes some labor to correctly define the multiplicity at 0 and infinity in the case of $(B_1, B_2)$, see definition \ref{defmultipl1}.
\end{rem}

This algebraic hypersurface is actively studied starting from the paper \cite{S94} and the book \cite{GKZ94}.  We aim at describing its singular locus $\sing R_B$.

In the classical case $B_i=\{0,1,2,\ldots,d_i\}$ the sparse resultant $R_B$ is the usual resultant, that is the hypersurface defined by Sylvester matrix. If $d_i > 2$, the singular locus has one irreducible component, and its codimension is 2. The irreducible component is % the space denoted below as $\widehat M_0^0(1,1)$, which for good $(B_1, B_2)$ % TODO: правильно? вроде кроме условий (3), (4)
% can be defined as
the closure of a (locally closed) irreducible set
$$H(1,1)=\{(f_1,f_2) \in\C^{B_1} \times \C^{B_2}\,|\,f_1(x)=f_2(x)=0\mbox{ has at least two non-zero solutions}\}.$$
Moreover, in the classical case we know how singular the resultant is at a generic point of its singular locus: its transversal singularity type at a generic point $(f_1,f_2)\in \sing R_B$ (i.e. the type of the singularity of the intersection of $R_B$ with a germ of a 2-dimensional plane transversal to $\sing R_B$ at $(f_1,f_2)$) is $\mathcal{A}_1$ (i.e. that of the union of two transversal lines in the plane).

Our main result describes the conditions on $B_1$ and $B_2$ under which similar description holds true. 

\begin{theor}\label{theormain}
i) There exists a codimension 3 subset $\Sigma\subset \C^{B_1}\x\C^{B_2}$ such that at every point of $\sing R_B\setminus\Sigma$ the transversal singularity type of $R_B$ is $\mathcal{A}_1$, unless at least one of the following conditions holds:

\begin{enumerate}
    \item One can shift $B_1$ and $B_2$ to the same proper sublattice $k\Z\subset\Z$ with $k\geqslant 2$.
    \item One can split one of $B_i$'s (say, $B_1$) into $B'\sqcup B''$ so that $B'$, $B''$ and $B_2$ can be shifted to the same sublattice $k\Z\subset\Z$ with $k\geqslant 3$.
    \item For every $i=1,2$, the two leftmost elements of $B_i$ differ by more than 1.
    \item For every $i=1,2$, the two rightmost elements of $B_i$ differ by more than 1.
    \item The number of elements in one of $B_i$'s (say, $B_1$) is 2, and $\max B_2-\min B_2 > 2$.
\end{enumerate}

ii) The singular locus of the resultant $R_B$ consists of several irreducible components of codimension 2, unless the following condition holds:

\begin{enumerate}
    \item[(6)] $B_1=\{i,i+k\}$ and $B_2=\{j,j+k\}$ for some $i,j \in Z$ and $k \in \N$.
\end{enumerate}
\end{theor}

% TODO: can sing R_B have more then one component if these conditions are not satisfied?

In fact, the condition (1) is redundant, because it implies the conditions (3) and (4), but it is convenient to the narrative. The following examples illustrate what happens to the singular locus once one of the conditions (1-6) takes place.
\begin{exa}\label{exaintro}
First let us consider part (i).

\medskip
1) If condition (1) takes place, then for every pair $(f_1, f_2)$ of Laurent polynomials with a non-zero common root $x$ there is also $(k-1)$ more non-zero common roots $x \cdot \epsilon_k$, where $\epsilon_k$ are $k$-th roots of unity.
In this case, the study of the resultant $R_B$ can be reduced to the study of a resultant $R_{B'}$ for a smaller pair of support sets $B_1'$ and $B_2'$ such that $B_1=k \cdot B_1'+m_1=\{kb+m: b \in B_1'\}$ and $B_2=k \cdot B_2'+m_2$. Indeed, $\C^{B_1} \times \C^{B_2} \cong \C^{B_1'} \times \C^{B_2'}$, and, moreover, $N_B(1)=N_{B'}(1)$, thus $R_B \cong R_{B'}$.

\medskip
2) If condition (2) takes place, then, for some component of the singular locus of the resultant, the transversal singularity type at a generic point is expected to have $k \geqslant 3$ components, thus to differ from $\mathcal{A}_1$.

For instance, consider $B_1=\{0,1,3\}$ and $B_2=\{0,3\}$. Let us denote the polynomials by $ax^3+bx+c$ and $dx^3+e$. Then the resultant is given by the equation
\begin{equation}\label{eqae_cd}
(ae-cd)^3+b^3d^2e=0
\end{equation}
This equation is homogeneous in $(a,b,c)$ and in $(d,e)$, so let us restrict it to $a=d=1$. Then the singular locus of the restriction is given by the equation $b=c-e=0$, which corresponds to a component of the singular locus of the subset defined by Formula (\ref{eqae_cd}). Choosing $b=\epsilon_1,\; c=1+\epsilon_2,\; e=1$ as a transversal plane germ to this component, we see that the transversal singularity of this component is given by $\epsilon_1^3-\epsilon_2^3=0$, thus it is the union of three transversal lines intersecting in a common point.

\medskip
3) If condition (3) or (4) takes place, then, for some component of the singular locus of the resultant, the transversal singularity type at a generic point is expected to have a singular component, thus to differ from $\mathcal{A}_1$.

For instance, consider the same example $B_1=\{0,1,3\}$ and $B_2=\{0,3\}$. There are two more components of the singularity locus which are given by $a=d=0$ and $d=e=0$. For the first one, the transversal singularity is the cusp $x^2=y^3$.

The case when condition (4) takes place can be reduced to (3) by the change of coordinate $\tilde x=x^{-1}$, because the sets $B_1$ and $B_2$ can be replaced by the sets $B_1'=-B_1=\{-b: b \in B_1\}$ and $B_2'=-B_2$.

\medskip
4) If condition (5) takes place, then again, for some component of the singular locus of the resultant, the transversal singularity type at a generic point is expected to have $k \geqslant 3$ components, thus to differ from $\mathcal{A}_1$.

For instance, suppose that $|B_1|=2$ and $B_2=\{0,1,\ldots,d_2\}$ with $d_2>2$. There is one more component of the singularity locus which is given by the equation $f \equiv 0$. This component has codimension 2, and in its general point one of the polynomials is identically 0 and the second has $d_2$ different roots, thus the transversal singularity type has $d_2>2$ components. % TODO: объединение гладких трансверсальных прямых?

The same happens with the component $d=e=0$ from the previous case.

\medskip
Now let us consider part (ii).

5) If condition (6) takes place, then the singular locus consists of one point $(0,0) \in \C^{B_1} \x\C^{B_2}$, thus has codimension 4. Indeed, without loss of generality we can replace $B_1$ and $B_2$ by $B'_1=B'_2=\{0,k\}$. The singular locus for $B'_1$ and $B'_2$ is the same as for $B''_1=B''_2=\{0,1\}$, because we consider them as subsets, not counting the multiplicities. Now let us denote the polynomials by $ax+b$ and $cx+d$. The resultant is given by the equation $ad=bc$, thus it has just one singular point $(0,0)$.

\end{exa}

% \begin{rem}
% As explained above, we can assume with no loss of generality that $B_1$ and $B_2$ cannot be shifted to a proper sublattice of $\Z$. We now assume this once and for all. % TODO: это используется в тексте ниже?
% \end{rem}

\begin{rem}
Actually, once the condition (2), (3), (4) or (5) takes place in the setting of the theorem, the singular locus of the resultant is always expected to have a component, at whose generic point the transversal singularity type differs from $\mathcal{A}_1$. However, the study of its singularity type (and the proof that it differs from $\mathcal{A}_1$) is non-trivial and will be done in a separate paper.
\end{rem}

\begin{rem}
While conditions similar to (1), (3) and (4) are familiar to the experts in tropical geometry and Newton polytopes, the condition (2) is new (to the best of our knowledge). The source of this condition is the subsequent theorem \ref{theorstratumranks21for111} describing degenerations of certain matrices of Vandermonde type.
\end{rem}
\begin{rem}
It would be important to obtain a version of this theorem for $A$-discriminants: their singular locus is studied e.g. in \cite{E13}, \cite{DHT16} and \cite{V21} under the assumption that it has the expected dimension and transversal singularity type, but so far there are no known criteria for these assumption to hold true even for univariate polynomials (except for some sufficient conditions in Section 3.4 of \cite{E13}).
\end{rem}
\begin{rem}
In this article, we prove Theorem \ref{theormain} over the field of complex numbers. This setting is important in our proof of Lemma \ref{lem3minor}, where one of the key steps essentially makes appeal to elementary geometry in the complex line (see Figure \ref{figinscrangles}). We shall address the case of arbitrary field in a subsequent paper, but we do not know to what extent it remains valid in the finite characteristics. % TODO: it is probably not important; what is probably important is that we speak about diffeomophisms
\end{rem}

% TODO: extend

Most of the paper is devoted to the proof of Theorem \ref{theormain}. The proof consists of two parts: in Sections \ref{secfiltr}--\ref{secstratrootspace}, we systematically study the natural stratifications of the resultant and certain related objects. The aim of this study is to reduce Theorem \ref{theormain} to several nontrivial facts about ranks of certain Vandermonde-type matrices. These facts (Lemma \ref{lemstratumranks11} and especially Theorem \ref{theorstratumranks21for111}) are then proved in Sections \ref{secstrats11forn111}--\ref{secstrats12forn111}.

In the last section, as an application of Theorem \ref{theormain}, we prove the following fact, partially answering the question given in the Remark 1.2 of \cite{V19}.

Given two finite sets $A_1$ and $A_2$ in $\Z^3$ and a pair of generic complex (or real) Laurent polynomials $f_1$ and $f_2$ supported at these sets, the equations $f_1=f_2=0$ define a smooth spatial algebraic curve in $\CC^3$ (or $(\R^*)^3$). The closure $C$ of its projection to the first coordinate plane is not in general smooth: at least it may have singularities of type $\mathcal{A}_1$ at the points having two preimages. Such $\mathcal{A}_1$ singularities are stable under local perturbations of the smooth spatial curve, similarly to the transversal self-intersections of a knot diagram.
\begin{theor}\label{theorproj} The curve $C$ has no other singularity types, unless the projections $B_1$ and $B_2$ of $A_1$ and $A_2$ to the last coordinate line satisfy one of the five conditions of the part i) of the main theorem \ref{theormain}.
\end{theor}

\begin{rem}
Even if the sets $A_1$ and $A_2$ in this theorem are the sets of all lattice points in prescribed Newton polytopes, their projections $B_1$ and $B_2$ rarely consist of several consequtive integers. They usually have gaps. This explains why we need Theorem 1.2 for arbitrary support sets, and not just for the classical case $B_i=\{0,1,2,\ldots,d_i\}$.
\end{rem}

\section{A filtration of the sparse resultant}\label{secfiltr}

This section is devoted to the first reductions of part (i) of the main theorem, as a result we reduce it to the much more concrete Theorem \ref{theorcodimNpre}, which we prove in the subsequent sections.

\subsection{A branched covering of the sparse resultant}

Let us recall that the sparse resultant $R_B\subset\C^{B_1} \times \C^{B_2}$ is the closure of the set
$$H(1)=\{(f_1,f_2)\in\C^{B_1} \times \C^{B_2}\,|\,f_1(x)=f_2(x)=0\mbox{ has at least one non-zero solution}\}.$$

Let us define the subset $\widetilde R_B \subset \C^{B_1}\x\C^{B_2}\x\CP^1$ to be the closure of
$$\widetilde H(1) = \{(f,g,x) \in \C^{B_1} \x \C^{B_2} \x \C^\x \,|\, f(x)=g(x)=0\}.$$
We would like to study $R_B$ using $\widetilde R_B$ and the tautological projection $\C^{B_1}\x\C^{B_2}\x\CP^1 \to \C^{B_1}\x\C^{B_2}$.

The subset $\widetilde H(1)$ is given by the equations $f(x)=\sum_{i \in B_1} f_i x^i=0$ and $g(x)=\sum_{j \in B_2} g_j x^j=0$. They are equations in variables $f_i$ for $i \in B_1$, variables $g_j$ for $j \in B_2$ and a variable $x$. Its closure $\widetilde R_B$ is the projectivization with respect to $x$, thus it is given by the equations
\begin{equation}\label{eqsheaf}
\widehat f(x,y)=\sum_{i \in B_1} f_i x^{i-\min B_1}y^{\max B_1-i}=0 \quad \text{ and } \quad \widehat g(x,y)=\sum_{i \in B_2} g_j x^{j-\min B_2}y^{\max B_2-j}=0,
\end{equation}
where $(x:y)$ are homogeneous coordinates on $\CP^1$. % TODO: (x:y) vs. x

\begin{lemma}\label{lemwidehatRBsmooth}
The set $\widetilde R_B$ is smooth. % and its intersection with $\C^{B_1} \times \C^{B_2} \times \C^\x$ is exactly $\widetilde H(1)$
\end{lemma}

\begin{proof}
{\ifshort
It follows from the fact that the space $\widetilde R_B$ is a vector bundle over $\CP^1$.
\fi}

{\iflong
Let us consider the subset $\C^{B_1}\x\C^{B_2}\x\C_y$, where $\C_y \subset \CP^1$ is defined by $y \neq 0$. As $\C^{B_1}\x\C^{B_2}\x\C_y \cong \C^{B_1}\x\C^{B_2}\x\C$, we have a canonical coordinate chart on it. For the subset $\C^{B_1}\x\C^{B_2}\x\C^\x$, where $\C^\x \subset \CP^1$ is defined by $x,y \neq 0$, it is true that the subset $\widetilde R_B \cap (\C^{B_1}\x\C^{B_2}\x\C^\x)=\widetilde H(1)$ is defined exactly by the equations $f(x)=g(x)=0$.

Let us suppose that $\min B_1=\min B_2=0$. Then for the canonical coordinate chart on $\C^{B_1}\x\C^{B_2}\x\C_y$ it is true that the subset $\widetilde R_B \cap (\C^{B_1}\x\C^{B_2}\x\C_y)$ is defined exactly by the equations $f(x)=g(x)=0$.

The differentials in this canonical coordinate chart are
$$df(x)=\sum_{i \in B_1} x^i df_i + \sum_{j \in B_2} 0\,dg_j + f'(x) dx,$$
$$dg(x)=\sum_{i \in B_1} 0\,df_i + \sum_{j \in B_2} x^j dg_j + g'(x) dx.$$

We supposed that $\min B_1=0$ and $\min B_2=0$, thus the differentials contain the summands $df_0$ and $dg_0$ and thus they are linearly independent at each point $(f,g,x) \in \C^{B_1}\x\C^{B_2}\x\C_y$, so that $\widetilde R_B \cap (\C^{B_1}\x\C^{B_2}\x\C_y)$ is smooth.

In general case, we can shift $B_1$ and $B_2$ inside $\Z$ and obtain the case when $\min B_1=\min B_2=0$. 

To prove the smoothness of $\widetilde R_B \cap (\C^{B_1}\x\C^{B_2}\x\C_x)$, where $\C_x \subset \CP^1$ is defined by $x \neq 0$, let us notice that it is equal to $\widetilde R_{B'} \cap (\C^{B_1'}\x\C^{B_2'}\x\C_y)$ for $B_1'=-B_1=\{-b_1 \,|\, b_1\in B_1\}$ and $B_2'=-B_2$. Consequently, the whole $\widetilde R_B \subset \C^{B_1}\x\C^{B_2}\x\CP^1$ is smooth.
\fi}
\end{proof}

\begin{lemma}
The subset $R_B$ is the image of $\widetilde R_B$ under the tautological projection $\C^{B_1}\x\C^{B_2}\x\CP^1 \to \C^{B_1}\x\C^{B_2}$.
\end{lemma}

\begin{proof}
The tautological projection maps $\widetilde H(1)$ to $H(1)$, thus it maps $\widetilde R_B$, which is the closure of $\widetilde H(1)$, to $R_B$, which is the closure of $H(1)$. By the definition of $\widetilde H(1)$ the restriction $\pi: \widetilde R_B \to R_B$ maps onto $H(1)$; the tautological projection is proper, thus the image of $\widetilde R_B$ is closed, thus the image is the whole $R_B$; that is, $\pi$ is surjective and $R_B$ is the image of $\widetilde R_B$.
\end{proof}

We will denote the restriction of the tautological projection as $\pi: \widetilde R_B \to R_B$.

As the set $\widetilde R_B$ is given by Formula (\ref{eqsheaf}), for each $(f,g) \neq (0,0)$ there is a finite number of preimages $(f,g,x)$. Consequently, $R_B$ can be non-smooth at the point $(f,g) \neq (0,0)$ in one of the following two cases: either there is more than one preimage $(f,g,x)$, or the map $\pi$ is not a local embedding at $(f,g,x)$. To detail it, we will need some definitions.

\subsection{Multiplicities of the roots}

Let $B \subset \Z$ be a finite set. Let $f(x)=\sum_{b\in B}c_bx^b \in\C^B$ be a (Laurent) polynomial with the support $B$ (below we will usually write simply "polynomial" instead of "Laurent polynomial" for conciseness). From now on, we shall identify $f(x)$ with the section
$$s_f=\sum_{b\in B}c_bx^{b-\min B}y^{\max B-b}$$
of the invertible sheaf $\mathcal{O}(\max B-\min B)$ on $\CP^1$ with homogeneous coordinates $x$ and $y$. % This identification depends on the choice of $B$.
Notice that the formula is the same as Formula (\ref{eqsheaf}).

The purpose is to be able to speak of the multiplicities of $f$ at $0$ and $\infty\in\CP^1$.

\begin{defin}\label{defmultipl1}
Let $f(x) \in \C^B$ be a polynomial, and $x \in \CP^1$ be a point. The multiplicity of $f$ at $x$ is the multiplicity of the root of the section $s_f$ at $x$. It is denoted $\ord_x(f)$. % or, if it is necessary to indicate $B$, $\ord_x^B(f)$.
\end{defin}

{\iflong
Explicitly, if $B=\{b_1, b_2, \ldots, b_m\}$ with $b_1 < b_2 < \ldots < b_m$, the multiplicity of $f$ at $x \in \C^\times$ is equal to the multiplicity of the root $x$ of the Laurent polynomial $f(x)$, and to calculate the multiplicity of $f$ at $0$ or $\infty$ we write down $c_{b_1}, 0, \ldots, 0, c_{b_2}, 0, \ldots, 0, c_{b_m}$, where $c_{b_i}$ is written on $(b_i-\min B_i+1)$-th position, and count the number of zeroes on the left or on the right respectively. If $f$ is identically zero, then the multiplicity of $f$ at any $x \in \CP^1$ is infinite. % TODO: правильно?
\fi}

\begin{rem}
If $f$ is a polynomial, not just a Laurent polynomial, then under this identification the multiplicity of the section $s_f$ at 0 is different from the multiplicity of the root of $f$ at 0 when $f$ is considered as a polynomial. {\iflong In particular, if we shift $B_i$ to the right by $n$ and multiply $f$ by $x^n$, then the multiplicity of $s_f$ at $0$ would not change, but the multiplicity of $f$ at 0 when $f$ is considered as a polynomial would increase by $n$.\fi}
\end{rem}

{\iflong
\begin{rem}
To study the multiplicity at 0, we could just shift $B$ such that $\min B=0$. But we also need the multiplicity at infinity, and the definition above allows us to do study both and in a systematic way.
\end{rem}
\fi}

\begin{defin}
Let $f_1(x) \in \C^{B_1}$ and $f_2(c) \in \C^{B_2}$ be two polynomials, and $x \in \CP^1$ be a point. The multiplicity of a pair $(f_1, f_2)$ at $x$ is the minimum of the multiplicities of the sections $s_{f_1}$ and $s_{f_2}$ at $x$. It is denoted $\ord_x(f_1, f_2)$. % or, if it is necessary to indicate $B_1$ and $B_2$, $\ord_x^{B_1,B_2}(f_1,f_2)$.
\end{defin}

In particular, we say that $(f_1,f_2)$ has a common root at $x \in \CP^1$ if $\ord_x(f_1,f_2) \geqslant 1$.

\begin{rem}
{\iflong For $x \in \C^\x$ it coincides with the usual notion of the common root of a pair of Laurent polymonials $f_1$ and $f_2$.\fi} For $x=0$ ($x=\infty$) it means that both $f_1$ and $f_2$ have zero leftmost (rightmost) coefficients respectively.
\end{rem}

% TODO: "root" or "common root"?

\subsection{Singular locus of the sparse resultant}\label{subsetsingloc}

The definition of common root of $(f,g)$ in $\CP^1$ allows us to see the set $\widetilde R_B$ as the set of triples $(f,g,x)$ where $x$ is a common root of $(f,g)$. As the sparse resultant $R_B$ is the image of $\widetilde R_B$, it consists of pairs $(f,g)$ that has a common root. This definition, more simple then the definition as the closure of $H(1)$, it possible due to introduction of the notion of common roots at 0 and infinity.

Now we also to see the map $\pi: \widetilde R_B \to R_B$ differently: the preimage $\pi^{-1}(f,g)$ consists of all $(f,g,x)$ such that $x$ is a common root of $(f,g)$. Thus if $(f,g)$ has $k$ common roots, then $R_B$ locally at $(f,g)$ consists of $k$ (or less) local branches corresponding to these roots. % TODO: why branches do not coincide?

\begin{lemma}\label{lemrbnonlocdiff}
The map $\widetilde R_B \to R_B$ is not local embedding at $(f,g,x)$ if and only if $x$ is a common root of $(f,g)$ of multiplicity more than one.
\end{lemma}

\begin{proof}
{\ifshort
It is a special case of Thom's transversality lemma. % TODO: details
\fi}

{\iflong
The map $\pi: \widetilde R_B \to R_B$ is local embedding at the point $(f,g,x)$ if and only if the tangent map $\pi_*: T_{(f,g,x)} \widetilde R_B \to T_{(f,g)}(\C^{B_1} \times \C^{B_2})$ is injective.

Suppose that $x \neq \infty$. Let us as in Lemma \ref{lemwidehatRBsmooth} consider the subset $\C^{B_1}\x\C^{B_2}\x\C_y$, where $\C_y \subset \CP^1$ is defined by $y \neq 0$, and the canonical coordinate chart on it. The tangent map of the projection $\C^{B_1} \times \C^{B_2} \times \C \to \C^{B_1} \times \C^{B_2}$ is the quotient by $\partial/\partial x$, and the tangent space $T_{(f,g,x)} \widetilde R_B$ is $V(df(x), dg(x))$, thus the tangent map $\pi_*$ is non-injective at $(f,g,x)$ if and only if $\partial/\partial x$ is the root of both $df(x)$ and $dg(x)$ at $(f,g,x)$.

From the explicit formulas for $df(x)$ and $dg(x)$ we see that it means $f'(x)=g'(x)=0$. The point $(f,g,x)$ is in $\widetilde R_B$, thus $f(x)=g(x)=0$, thus the map $\pi: \widetilde R_B \to R_B$ is not a local embedding at $(f,g,x)$ if and only if $x$ is a root of $(f, g)$ of multiplicity at least two.

If $x = \infty$, then we can replace $B_i$ with $B_i'=-B_i$ and prove the same.
\fi}
\end{proof}

\begin{sledst}\label{corsing}
Suppose that $B_1$ and $B_2$ can not be shifted to the same proper sublattice $k\Z\subset\Z$ with $k\geqslant 2$. The subvariety $R_B$ is singular at the point $(f,g)\in \C^{B_1}\x\C^{B_2}$ % different from $(0,0)$
if and only if

1) $(f,g)$ has a root in $\CP^1$ of multiplicity at least 2, or

2) $(f,g)$ has at least two roots in $\CP^1$.
\end{sledst}

\begin{proof}
Indeed, by Lemma \ref{lemwidehatRBsmooth} $R_B$ is the image of a smooth set $\widetilde R_B$. The set $R_B$ can be non-smooth at the point $(f,g)$ in one of the following two cases: either the map $\pi$ is not a local embedding at $(f,g,x)$, or there is more than one preimage $(f,g,x)$. It corresponds to a root of multiplicity at least 2 and two roots respectively.

Vice versa, suppose that $(f,g) \neq (0,0)$ has a root of multiplicity at least 2. Then the map $\pi$ is not a local embedding at $(f,g,x)$. This map has finite fibers (expect for $(0,0)$), thus the image $R_B$ is singular at $(f,g)$.

Now suppose that $(f,g) \neq (0,0)$ has at least two roots. If some of them have multiplicity at least 2, then $R_B$ is singular at $(f,g)$. If all have multiplicity 1, then $R_B$ at the point $(f,g)$ has several local branches corresponding to the roots, each of which is smooth. They can not coincide. Indeed, if two branches coindice in a small neighborhood $U$ of $(f,g)$ in $R_B$, then every pair $(f_1,g_1)\in U$ has two common roots, which is possible only if $B_1$ and $B_2$ can be shifted to the same proper sublattice $k\Z\subset\Z$ with $k\geqslant 2$.
% TODO: more details

At last, if $(f,g)=(0,0)$, then then the conical space $R_B$ can be smooth if and only if it is linear, that can not happen. {\iflong
Indeed, for a linear set the equations $df(x)$ and $dg(x)$ that define the tangent space to $R_B$ should be constant, which means $|B_1|=|B_2|=1$.
\fi}
\end{proof}

\begin{lemma}\label{lemtransvers}
Suppose that $B_1$ and $B_2$ can not be shifted to the same proper sublattice $k\Z\subset\Z$ with $k\geqslant 2$. Let $x_1$ and $x_2$ be different roots of $(f,g)$ of multiplicity 1. Then the corresponding local branches of $R_B$ are transversal hypersurfaces, unless $x_1$ and $x_2$ are both from $\C^\x$ and either both have multiplicity at least 2 for $f$, or both have multiplicity at least 2 for $g$.
\end{lemma}

\begin{proof}
Suppose that neither $x_1$, nor $x_2$ is equal to the infinity. {\ifshort
Let us shift $B_1$ and $B_2$ so that $\min B_1=\min B_2=0$, and consider the subset $\C^{B_1}\x\C^{B_2}\x\C_y$, where $\C_y \subset \CP^1$ is defined by $y \neq 0$. As $\C^{B_1}\x\C^{B_2}\x\C_y \cong \C^{B_1}\x\C^{B_2}\x\C$, we have a canonical coordinate chart on it. In this coordinate chart, the subset $\widetilde R_B \cap (\C^{B_1}\x\C^{B_2}\x\C_y)$ is defined exactly by the equations $f(x)=g(x)=0$.
\fi}
{\iflong
Let us as in Lemma \ref{lemwidehatRBsmooth} shift $B_1$ and $B_2$ so that $\min B_1=\min B_2=0$, and consider the subset $\C^{B_1}\x\C^{B_2}\x\C_y$, where $\C_y \subset \CP^1$ is defined by $y \neq 0$, and the canonical coordinate chart on it.
\fi}

{\ifshort
The tangent space to $R_B$ is given by the equations
$$df(x)=\sum_{i \in B_1} x^i df_i + \sum_{j \in B_2} 0\,dg_j + f'(x) dx,$$
$$dg(x)=\sum_{i \in B_1} 0\,df_i + \sum_{j \in B_2} x^j dg_j + g'(x) dx.$$

So now consider the map $\pi: \widetilde R_B \to R_B$ and its tangent map $\pi_*: V(df(x), dg(x)) \to  T_{(f,g)}(\C^{B_1} \times \C^{B_2})$ at the points $(f,g,x_k)$, $k=1,2$.
\fi}
{\iflong
Similarly to Lemma \ref{lemrbnonlocdiff}, consider the map $\pi: \widetilde R_B \to R_B$ and its tangent map $\pi_*: V(df(x), dg(x)) \to  T_{(f,g)}(\C^{B_1} \times \C^{B_2})$ at the points $(f,g,x_k)$, $k=1,2$.
\fi}
We have that either $f'(x_k) \neq 0$ or $g'(x_k) \neq 0$, thus the image is given by the linear equation
$$g'(x_k) \cdot df(x)-f'(x_k) \cdot dg(x)=g'(x_k) \cdot \sum_{i \in B_1} x_k^i df_i - f'(x_k) \cdot \sum_{j \in B_2} x_k^j dg_j.$$

Suppose that the hypersurfaces are not transversal, thus such covectors for $x_1$ and $x_2$ are proportional and non-zero. The covectors have summands $f'(x_1)x_1^mdg_m$ and $f'(x_2)x_2^mdg_m$ with $m > 0$, thus if $x_1=0$ and $x_2 \neq 0$, then $f'(x_2)=0$, and similarly $g'(x_2)=0$. The point $(f,g,x_2)$ is in $\widetilde R_B$, thus $f(x_2)=g(x_2)=0$ and so $x_2$ is a root of $(f,g)$ of multiplicity at least 2, but that contradicts the conditions of the lemma. Thus $x_1=0$ and $x_2 \neq 0$ can not happen, and as $x_1 \neq x_2$ and the conditions are simmetric in $x_1$ and $x_2$, we have that $x_1 \neq 0$ and $x_2 \neq 0$.

The covectors have summands $f'(x_1)dg_0$ and $f'(x_2)dg_0$, thus either both $f'(x_1)$ and $f'(x_2)$ vanish, or both are non-zero. The points $(f,g,x_k)$ are in $\widetilde R_B$, thus $f(x_1)=f(x_2)=0$. By the conditions of the lemma, neither $x_1$ and $x_2$ have multiplicity at least 2 for $f$, thus $f'(x_1) \neq 0$ and $f'(x_2) \neq 0$. Consequently, the covectors $\sum_{j \in B_2} x_1^j dg_j$ and $\sum_{j \in B_2} x_2^j dg_j$ are proportional, thus $x_1/x_2$ is a $k$-th root of unity with $k\geqslant 2$ and one can shift $B_2$ to a proper sublattice $k\Z\subset\Z$. Similarly with $g'$ and $B_1$, thus $B_1$ and $B_2$ can be shifted to the same proper sublattice $k\Z\subset\Z$ with $k\geqslant 2$, which contradicts the conditions of the lemma.

If neither $x_1$, nor $x_\infty$ is equal to 0, then we can replace $B_i$ with $B_i'=-B_i$ and prove the same.

If $x_1=0$ and $x_2=\infty$, then the equations of the hypersurfaces are
$$g'(x_1)df_{\min B_1}-f'(x_1)dg_{\min B_2} \text{ and } g'(x_2)df_{\max B_1}-f'(x_2)dg_{\max B_2},$$
thus they are also transversal. % TODO: details
\end{proof}

\begin{sledst}\label{cortranstype}
Suppose that $B_1$ and $B_2$ can not be shifted to the same proper sublattice $k\Z\subset\Z$ with $k\geqslant 2$. The transversal singularity type of $R_B$ at the singular point $(f,g) \in \C^{B_1}\x\C^{B_2}$ is $\mathcal A_1$ unless

1) $(f,g)$ has three roots in $\CP^1$,

2) $(f,g)$ has two roots in $\C^\x$, both of which have multiplicity at least two for $f$ (or, similarly, for $g$),

3) $(f,g)$ has a root of multiplicity at least 2.
\end{sledst}

In particular, it takes into consideration the case $(f,g)=(0,0)$.

% TODO: remark that it is one way, see N(22,11)

\subsection{The filtration}

The definition of the multiplicities allows us to filter the space $\C^{B_1} \times \C^{B_2}$, containing the resultant $R_B$, into the subsets $N(p)$ defined by the orders of common roots. It will be more convenient to narrow them down to
$$\C^{B_1\x} \x \C^{B_2\x}=\left\{ (f_1,f_2) \in \C^{B_1} \x \C^{B_2} \,|\, f_1 \not\equiv 0 \text{ and } f_2 \not\equiv 0 \right\},$$
because if $f_i$ is identically zero, then it has more roots than is expected. In particular, the subsets $N(1)$ and $N(1,1)$ defined below are the intersections of the subsets $H(1)$ and $H(1,1)$ from the introduction and the beginning of this section with the subset $\C^{B_1\x}\x\C^{B_2\x}$.

\begin{defin}\label{defN}
A symmetric filtration subset
$$N_{j_0}^{j_\infty}(j_1,\ldots,j_k),\quad j_0, j_\infty \geqslant 0, \quad j_1,\ldots,j_k \geqslant 1, \quad k \geqslant 0$$
consists of $(f_1,f_2) \in \C^{B_1\x} \times \C^{B_2\x}$ such that
\begin{enumerate}
    \item $\ord_0(f_1, f_2) \geqslant j_0$ and $\ord_\infty(f_1, f_2) \geqslant j_\infty$;
    \item $f_1=f_2=0$ has at least $k$ distinct solutions $x_1,\ldots,x_k$ in $\C^\times$;
    \item at the $m$-th solution $x_m$ holds $\ord_{x_m}(f_1, f_2) \geqslant j_m$.
\end{enumerate}
\end{defin}

The numbering of the symmetric filtration subsets is defined up to a permutation of $(j_1, \ldots, j_k)$: for example, $N_1^2(3,4,5)$ is equal to $N_1^2(4,5,3)$, but in general not to $N_2^1(3,4,5)$.

If $k=0$, the symmetric filtration subset is denoted simply $N_{j_0}^{j_\infty}$. We will also denote $N_0^0(j_1,\ldots,j_k)$ by $N(j_1,\ldots,j_k)$.

The subset $N(1)$ is $H(1) \cap (\C^{B_1\x}\x\C^{B_2\x})$, thus the closure of $N(1)$ is $R_B \cap (\C^{B_1\x}\x\C^{B_2\x})$.

Now we can reformulate the $\C^{B_1\x}\x\C^{B_2\x}$-part of Corollary \ref{corsing} as follows:

\begin{sledst}
Suppose that $B_1$ and $B_2$ can not be shifted to the same proper sublattice $k\Z\subset\Z$ with $k\geqslant 2$. Then the intersection of the singular locus $\sing R_B$ with $\C^{B_1\x}\x\C^{B_2\x}$ is
$$N(2) \cup N^2_0 \cup N^0_2 \cup N(1,1) \cup N^1_0(1) \cup N^0_1(1).$$
\end{sledst}

Moreover, we can generalize the last definition for the case of two separate lines of indices.

\begin{defin}\label{defNN}
A general filtration subset
$$N\;_{j_0}^{j_\infty}\left(\,^{j^1_1,\ldots,j^1_k}_{j^2_1,\ldots,j^2_k}\right),\quad j_0,j_\infty \geqslant 0, \quad j^1_1,\ldots,j^1_k,j^2_1,\ldots,j^2_k \geqslant 1$$
consists of $(f_1,f_2) \in \C^{B_1\x} \times \C^{B_2\x}$ such that
\begin{enumerate}
    \item $\ord_0(f_1, f_2) \geqslant j_0$ and $\ord_\infty(f_1, f_2) \geqslant j_\infty$;
    \item $f_1=f_2=0$ has at least $k$ distinct solutions $x_1,\ldots,x_k$ in $\C^\times$;
    \item at the $m$-th solution $x_m$ holds $\ord_{x_m}(f_1) \geqslant j^1_m$ and $\ord_{x_m}(f_2) \geqslant j^2_m$.
\end{enumerate}

The numbering of the general filtration subsets is defined up to a simultaneous permutation of $(j^1_1, \ldots, j^1_k)$ and $(j^2_1, \ldots, j^2_k)$.

We will denote $N\;_0^0\left(\,^{j^1_1,\ldots,j^1_k}_{j^2_1,\ldots,j^2_k}\right)$ by $N\left(\,^{j^1_1,\ldots,j^1_k}_{j^2_1,\ldots,j^2_k}\right)$.
\end{defin}

Here we have that
$$N_{j_0}^{j_\infty}(j_1,\ldots,j_k)=N\;_{j_0}^{j_\infty}\left(\,^{j_1,\ldots,j_k}_{j_1,\ldots,j_k}\right).$$

% TODO: N not locally closed

Now we can similarly reformulate the $\C^{B_1\x}\x\C^{B_2\x}$-part of Corollary \ref{cortranstype}:
\begin{sledst}\label{cortranstype2}
The transversal singularity type of $R_B$ at the point $(f,g) \in \C^{B_1\x}\x\C^{B_2\x}$ is $\mathcal A_1$ unless $(f,g)$ lies in one of the following filtration subsets:

1) $N(1,1,1)$, $N^0_1(1,1)$, $N^1_0(1,1)$ or $N_1^1(1)$;

2) $N(2)$, $N_2^0$ or $N_0^2$;

3) $N\left(\,^{2,2}_{1,1}\right)$ or $N\left(\,^{1,1}_{2,2}\right)$
\end{sledst}

\begin{rem}\label{remideaofmainthm}
One can easily see that almost all of that subsets $N(p)$ listed above usually have codimension at least 3. Moreover, we have that $N\left(\,^{2,2}_{1,1}\right) \subset N\left(\,^{2,1}_{1,1}\right)$ and $N\left(\,^{1,1}_{2,2}\right) \subset N\left(\,^{1,1}_{2,1}\right)$. Consequently, if the subsets $N(1,1,1)$, $N\left(\,^{2,1}_{1,1}\right)$ and $N\left(\,^{1,1}_{2,1}\right)$ have codimension at least 3, then there exists a codimension 3 subset $\Sigma\subset \C^{B_1}\x\C^{B_2}$ such that at every point of $(\sing R_B\setminus\Sigma) \cap (\C^{B_1\x}\x\C^{B_2\x})$ the transversal singularity type of $R_B$ is $\mathcal{A}_1$. It is the part (i) of the main theorem \ref{theormain} expect for the subset $(\{0\} \x \C^{B_2}) \cup (\C^{B_1} \x \{0\})$. % TODO: to be clear about codimensions
\end{rem}

\subsection{Expected codimensions of the filtration subsets}

The relation of the subsets $N(p)$ to the singular locus of the resultant $R_B$ motivates our interest in their codimension.

\begin{defin}
The expected codimension of a general fitration subset $N\;_{j_0}^{j_\infty}\left(\,^{j^1_1,\ldots,j^1_k}_{j^2_1,\ldots,j^2_k}\right)$ is the number
$$2j_0+2j_\infty+\sum_{m=1}^k (j^1_m+j^2_m-1).$$

In particular, the expected codimensions of a symmetric filtration subset $N_{j_0}^{j_\infty}(j_1,\ldots,j_k)$ is
$$2j_0+2j_\infty+\sum_{m=1}^{k} (2j_m-1).$$
\end{defin}

{\iflong
\begin{rem}
Indeed, to fix orders of the root at 0, one should impose $2j_0$ conditions on $(f_1, f_2)$: $j_0$ conditions that $f_1$ has a root of order at least $1, 2, \ldots, j_0$ at 0, and the same with $f_2$; similarly with $\infty$; and to fix a common root that has order $j^1_m$ for $f_1$ and $j^2_m$ for $f_2$ and is lying in $\C^\times$, one should first impose the condition of the existence of a common root $x_m$, then impose $(j^1_m-1)$ conditions that $f_1$ has a root of order at least $2, \ldots, j_m$ at $x_m$, and the same with $f_2$. We expect that each condition increases the codimension by 1.

% Moreover, one should impose the conditions that $f_1$ and $f_2$ have no other roots, and the condition that for a common root of $f_1$ and $f_2$ its order is not greater than expected for both $f_1$ and $f_2$ at once; but these conditions are expected to be open.
\end{rem}
\fi}

\begin{exa}
The symmetric filtration subsets of expected codimensions 1, 2 and 3 are
\begin{enumerate}
    \item $N(1)$;
    \item $N_1^0$, $N_0^1$ and $N(1,1)$;
    \item $N_1^0(1)$, $N_0^1(1)$, $N(2)$ and $N(1,1,1)$.
\end{enumerate}
\end{exa}

\begin{rem}
The filtration subsets % $N(1,1,1)$,
$N\left(\,^{2,1}_{1,1}\right)$ and $N\left(\,^{1,1}_{2,1}\right)$ from Remark \ref{remideaofmainthm} have expected codimension 3.
\end{rem}

\begin{rem}\label{remcodimreduc}
An algebraic subset can consist of several irreducible components of different dimensions, thus we should be accurate when talking about its codimension. When we say that an algebraic subset has codimension equal to the given number, we mean that the codimension of every its irreducible component is equal to the given number. When we say that the codimension of an algebraic subset is at most (at least) the given number, we mean that the codimension of every its irreducible component is at most (at least) the given number. In particular, in all of the mentioned cases it may be empty, although some sources consider the empty set to have codimension $\infty$.
\end{rem} % TODO: non locally closed - changes anything?

\subsection{Proof of the part (i) of the main theorem}

We have a following theorem:
\begin{theor}\label{theorcodimNpre}
Unless $(B_1,B_2)$ satisfy one of the five conditions of the part (i) of the main theorem \ref{theormain}, every symmetric filtration subset $N(p)$ of the expected codimension 1, 2 or 3 has actual codimension at least 1, 2 or 3 respectively.

Moreover, under these conditions, the filtration subsets $N\left(\,^{2,1}_{1,1}\right)$ and $N\left(\,^{1,1}_{2,1}\right)$ also have codimension at least 3.
\end{theor}

\begin{rem}
Each of the conditions of the part (i) of the main theorem \ref{theormain} is only used to estimate the codimensions of some of the strata. The exact details see in the theorem 
\ref{theorcodimN}.
\end{rem}

We will prove this theorem in Sections \ref{secstratrootspace}-\ref{secstrats12forn111}. Now we will use this theorem and Corollary \ref{cortranstype2} to prove the part (i) of the main theorem \ref{theormain}.

\begin{proof}[Proof of the main theorem \ref{theormain}, part (i)]
Let us first consider the $\C^{B_1\x}\x\C^{B_2\x}$-part. By Theorem \ref{theorcodimNpre}, in the conditions of the part (i) of the main theorem the filtration subsets
$$N(1,1,1), \quad N(2), \quad N\left(\,^{2,1}_{1,1}\right), \quad N\left(\,^{1,1}_{2,1}\right), \quad N_0^1(1) \; \text{ and } \; N_1^0(1)$$
have codimensions at least three, thus the same is true for their subsets
$$N(1,1,1), \quad N(2), \quad N\left(\,^{2,2}_{1,1}\right), \quad N\left(\,^{1,1}_{2,2}\right), \quad N_0^1(1,1) \; \text{ and } \; N_1^0(1,1).$$
The same is obviously true for the rest of the filtration subsets of Corollary \ref{cortranstype2}, namely the subsets 
$$N_1^1(1), \quad N_2^0 \; \text{ and } \; N_0^2,$$
thus by Corollary \ref{cortranstype2} there exists a codimension 3 subset $\Sigma'\subset \C^{B_1\x}\x\C^{B_2\x}$ such that at every point of $(\sing R_B \cap (\C^{B_1\x}\x\C^{B_2\x}))\setminus\Sigma'$ the transversal singularity type of $R_B$ is $\mathcal{A}_1$. % TODO: to be clear about codimensions

\medskip
Let us now consider the $(\{0\} \x \C^{B_2}) \cup (\C^{B_1} \x \{0\})$-part. We need to find a codimension 3 subset $\Sigma''\subset (\{0\} \x \C^{B_2}) \cup (\C^{B_1} \x \{0\})$ such that at every point of $(\sing R_B \cap (\C^{B_1}\x\C^{B_2})) \setminus\Sigma''$ the transversal singularity type of $R_B$ is $\mathcal{A}_1$.

By Corollary \ref{corsing}, the subset $\{0\}\x\C^{B_2\x}$ appears in $\sing R_B$ if and only if $B_2 \neq \{j, j+1\}$. Its codimension if $|B_1|$, thus if $|B_1|>2$, then we can include $(\{0\}\x\C^{B_2\x})$ into $\Sigma''$. Otherwise $|B_1|=2$ and it is a subset of codimension 2. By Lemma \ref{lemtransvers}, at its general point $(0,g)\in \{0\}\x\C^{B_2\x}$ (such that $g$ has no multiple roots) the subset $R_B$ has $\max B_2-\min B_2$ smooth transversal branches, thus if $\max B_2-\min B_2=2$, we can still add to $\Sigma''$ a part of $\{0\}\x\C^{B_2\x}$ containing polynomials with multiple roots. The case of $|B_1|=2$ and $\max B_2-\min B_2>2$ is excluded by the condition (5). % TODO: rewrite

The same can be proven about $\C^{B_1\x} \x \{0\}$.

As for the subset $\{(0,0)\}$, it has codimension at least 4.

Thus we get a subset $\Sigma''$ as necessary.

\medskip
We can now take $\Sigma$ to be $\Sigma' \cup \Sigma''$.
\end{proof}

\section{A stratification of the sparse resultant}

\subsection{The stratification}

In this section, we will study the strata $M(p)$ of a stratification related to the filtration $N$ introduced in Section \ref{secfiltr}. It will be more convenient again to narrow them down to
$$\C^{B_1\x} \x \C^{B_2\x}=\left\{ (f_1,f_2) \in \C^{B_1} \x \C^{B_2} \,|\, f_1 \not\equiv 0 \text{ and } f_2 \not\equiv 0 \right\},$$
because if $f_i$ is identically zero, then it has more roots than is expected.

\begin{defin}
A stratum
$$M_{j_0}^{j_\infty}(j_1,\ldots,j_k),\quad j_0, j_\infty \geqslant 0, \quad j_1,\ldots,j_k \geqslant 1, \quad k \geqslant 0$$
consists of $(f_1,f_2) \in \C^{B_1\x} \times \C^{B_2\x}$ such that
\begin{enumerate}
    \item $\ord_0(f_1, f_2) = j_0$ and $\ord_\infty(f_1, f_2) = j_\infty$;
    \item $f_1=f_2=0$ has exactly $k$ distinct solutions $x_1,\ldots,x_k$ in $\C^\times$;
    \item at the $m$-th solution $x_m$ holds $\ord_{x_m}(f_1, f_2) = j_m$.
\end{enumerate}
\end{defin}

The strata numbering is defined up to a permutation of $(j_1, \ldots, j_k)$: for example, $M_1^2(3,4,5)$ is equal to $M_1^2(4,5,3)$, but not to $M_2^1(3,4,5)$. The strata that do no differ by the permutation of $(j_1, \ldots, j_k)$ do not intersect. The space $\C^{B_1\x} \x \C^{B_2\x}$ is the disjoint union of such strata.

If $k=0$, a stratum is denoted simply as $M_{j_0}^{j_\infty}$.

% Our interest to these strata comes from their relation to the singular locus of the resultant $R_B$, or rather the part lying in $R_B \cap (\C^{B_1\x}\x\C^{B_2\x})$.

To study these sets, let us first notice that these sets are indeed strata in the sense that they are locally (Zariski) closed. To prove this, we introduce the following partial order $\succcurlyeq$ on the tuples $I=(j_0; j_\infty; j_1, \ldots, j_k)$ indexing these strata.
\begin{defin}\label{defstratorder}
Let us consider two strata
$$M(p)=M_{j_0}^{j_\infty}(j_1,\ldots,j_k) \quad \text{and} \quad M(q)=M_{g_0}^{g_\infty}(g_1,\ldots,g_l).$$
We say that $q \succcurlyeq p$ if, informally, the common roots of a pair $(h_1,h_2)\in M(q)$ in $\CP^1$ can be obtained from the common roots of a pair $(f_1,f_2)\in M(p)$ in $\CP^1$ by
\begin{enumerate}
    \item increasing their multiplicities,
    \item gluing some of the roots in $\C^\times$ together,
    \item moving some of the roots from $\C^\times$ to $0$ ot $\infty$ (but not moving the roots at $0$ or $\infty$ anywhere),
    \item adding some new common roots.
\end{enumerate}

Formally, $q \succcurlyeq p$ if for each $i \in \{0; \infty; 1, \ldots, k\}$ there is $r(i)\in \{0; \infty; 1, \ldots, l\}$ such that $r(0)=0$, $r(\infty)=\infty$ and for any $s \in \{0, \infty, 1, \ldots, l\}$ holds
$$\left(\sum_{i: \; r(i)=s} j_i\right) \leqslant g_s.$$
\end{defin}

\begin{defin}\label{defmwidehat}
Let us also define closed strata to be
$\widehat M(p)=\bigsqcup_{q \succcurlyeq p} M(q)$.
\end{defin}

\begin{exa}\label{exawidehat0011}
The closed stratum $\widehat M(1)$ consists of all pairs $(f_1, f_2)$ (of non-zero polynomials) that have at least one common root in $\CP^1$ of multiplicity at least 1. By the remark at the beginning of Subsection \ref{subsetsingloc}, $\widehat M(1)=R_B \cap (\C^{B_1\x}\x\C^{B_2\x})$. Moreover, in terms of Section \ref{secfiltr}, $\widehat M(1)=N(1)\cup N^1\cup N_1$. 

Similarly, the closed stratum $\widehat M(1,1)$ consists of all pairs $(f_1, f_2)$ that have either at least two common roots in $\CP^1$ or at least one common root in $\CP^1$ of multiplicity at least 2. By Corollary \ref{corsing}, if $B_1$ and $B_2$ can not be shifted to the same proper sublattice $k\Z\subset\Z$ with $k\geqslant 2$, then $\widehat M(1,1)=\sing R_B \cap (\C^{B_1\x}\x\C^{B_2\x})$. In terms of Section \ref{secfiltr}, $\widehat M(1,1)=N(2) \cup N^2 \cup N_2 \cup N(1,1) \cup N^1(1) \cup N_1(1)$.
\end{exa}

In Subsection \ref{subsecclosstrat} we will prove the following lemma, see Corollary \ref{corcodimwidehatM}:
\begin{lemma} The closed strata
$\widehat M(p)=\bigsqcup_{q \succcurlyeq p} M(q)$ are indeed (Zariski) closed in $\C^{B_1\x}\x\C^{B_2\x}$.
\end{lemma}

\begin{sledst}
The sets $M(p)=\widehat M(p)\setminus\bigcup_{q \succ p} \widehat M(q)$ are locally closed.
\end{sledst}

\begin{rem}\label{remb1ab}
The closed stratum $\widehat M(p)$ is not always equal to the closure of $M(p)$.

For example, let us consider $B_1=\{a, b\}$ and any $B_2$ (with at least two elements, as usual). A polynomial $f_1 \in \C^{B_1\x}$ has form $c_1x^a+c_2x^b$ and thus can not have two equal roots from $\C^\times$, thus the stratum $M_0^0(2)=\varnothing$. At the same time $\widehat M_0^0(2)$ contains $\widehat M_2^0$, which is non-empty if $b-a \geqslant 2$.

As we can see in this example, a stratum may be empty even if its expected codimension is less than the dimension of ambient space. Indeed, $M_0^0(2)$ has expected codimension 3 and lies in the space $\C^{B_1\x} \times \C^{B_2\x}$ of the dimension at least 4.
\end{rem}

\begin{rem}
The symmetric filtration subsets $N(p)$ are different from $M(p)$ in the fact that the conditions are inequalities, not equalities. But they are also different from $\widehat M(p)$, because, saying informally, inside $\widehat M(p)$ one can move points to 0 or $\infty$ and glue points. For example, $M_1^0(1)$ and $M(2)$ lie in $\widehat M_0^0(1,1)$, but usually not in $N_0^0(1,1)$.
\end{rem}

{\iflong
\begin{rem}
The subset $M(p)$ may be empty, while $N(p)$ is not empty. For example, if $B_1=kB_1'$ and $B_2=kB_2'$ with $k \geqslant 2$, then similarly to Example \ref{exaintro}, part 1), $M_0^0(1)=\varnothing$ and $N_0^0(1)=N_0^0(k)$ has codimension 1.

The subset $N(p)$ may be empty, while $\widehat M(p)$ is not empty. For example, in the conditions of Remark \ref{remb1ab} $N_0^0(2)=\varnothing$ and $\widehat M_0^0(2) \supset \widehat M_2^0$, which is non-empty if $b-a \geqslant 2$. In particular, $\widehat M(p)$ is not always equal to the closure of $N(p)$.
\end{rem}
\fi}

\subsection{Expected codimensions of the strata}

\begin{defin}
Just like for $N_{j_0}^{j_\infty}(j_1,\ldots,j_k)$, the expected codimension of the subset $M_{j_0}^{j_\infty}(j_1,\ldots,j_k)$ is the number
$$2j_0+2j_\infty+\sum_{m=1}^{k} (2j_m-1).$$
\end{defin}

\begin{rem}
If $q \succcurlyeq p$, then the expected codimension of $M(q)$ is not greater than the expected codimension of $M(p)$. If $q \succcurlyeq p$ and the expected codimensions or $M(q)$ and $M(p)$ are equal, then $q=p$.
\end{rem}

\begin{lemma}\label{lemcodimclassic}
In the classical case $B_i=\{0,1,2,\ldots,d_i\}$, the codimension of the stratum $M(p)$ (at every its irreducible component) is equal to the expected one.
\end{lemma}
For the sake of completeness, we provide the proof of this classical fact by parameterizing naturally a finite cover of the set $M(p)$.
\begin{proof}
% TODO: roots at infinity?

First, let us consider the case of $M(p)=M^0_0(j_1,\ldots,j_k)$. The idea is that, as neither $f_1$, nor $f_2$ is identically zero, the condition $(f_1, f_2) \in M^0_0(j_1,\ldots,j_k)$ is equivalent to the condition that they can be written as
$$f_1(t)=\prod_{i=1}^k (t-z_i)^{j_i} \prod_{i=1}^{e_1} (t-x_i) \text{ and } f_2(t)=\prod_{i=1}^k (t-z_i)^{j_i} \prod_{i=1}^{e_2} (t-y_i)$$
with some open conditions on $x_i$, $y_i$ and $z_i$, and that there is only a finite number of ways to write $(f_1, f_2)$ in this form. It works only for $B_i$ without gaps, because otherwise not all $(f_1, f_2)$ of the form above will lie in $\C^{B_1\x} \times \C^{B_2\x}$.

{\iflong
Formally, consider the subset
$$P \subset \C^{e_1} \times \C^{e_2} \times \C^k, \quad e_1=d_1-\sum_{m=1}^k j_m, \quad e_2=d_2-\sum_{m=1}^k j_m$$
of tuples $(x_1,\ldots,x_{e_1}; \; y_1,\ldots,y_{e_2}; \; z_1,\ldots,z_k)$ such that $z_i$ are pairwise different and non-zero, and that for each $m=1,\ldots,k$ at least one of the polynomials $\prod_{i=1}^{e_1}(x_i-z_m)$ and $\prod_{i=1}^{e_1}(x_i-z_m)$ is non-zero. It is an open subset of $\C^{d_1} \times \C^{e_2} \times \C^k$ and thus it has the dimension $d_1+d_2-\sum_{m=1}^{k} (2j_m-1)$.

There is a map $P \to M(p)$ of the form
$$\quad (x_1,\ldots,x_{e_1}; \; y_1,\ldots,y_{e_2}; \; z_1,\ldots,z_k) \mapsto \left( \prod_{i=1}^{e_1} (t-x_i) \prod_{i=1}^k (t-z_i)^{j_i}, \; \prod_{i=1}^{e_2} (t-y_i) \prod_{i=1}^k (t-z_i)^{j_i} \right),$$
which is surjective and has finite fibers, thus $\dim P=\dim M(p)$ and
$$\codim M(p)=d_1+d_2-\dim M(p)=\sum_{m=1}^{k} (2j_m-1),$$
which is the expected codimension.

To prove the general case, let us notice that the stratum
$$M_{j_0}^{j_\infty}(j_1,\ldots,j_k) \subset \C^{B_1\x} \times \C^{B_2\x} \text{ with } B_i=\{0,\ldots,d_i\}$$
is identified with the stratum
$$M_0^0(j_1,\ldots,j_k) \subset \C^{B_1'\x} \times \C^{B_2'\x} \text{ with } B_i'=\{j_0,\ldots,d_i-j_\infty\}$$
under the natural embedding $\C^{B_1'\x} \times \C^{B_2'\x} \hookrightarrow \C^{B_1\x} \times \C^{B_2\x}$, and that the last stratum can identified with the stratum
$$M_0^0(j_1,\ldots,j_k) \subset \C^{B_1''\x} \times \C^{B_2''\x} \text{ with } B_i''=\{0,\ldots,d_i-j_0-j_\infty\}$$
by the inclusion $\C^{B_1''\x} \times \C^{B_2''\x} \hookrightarrow \C^{B_1'\x} \times \C^{B_2'\x}$ induced by the pair of maps $B_i'' \to B_i'$, $b \mapsto b+j_0$.

Thus
\begin{equation*}
\begin{split}
&\codim M_{j_0}^{j_\infty}(j_1,\ldots,j_k)=d_1+d_2-\dim M_{j_0}^{j_\infty}(j_1,\ldots,j_k)=\\
&=d_1+d_2-((d_1-j_0-j_\infty)+(d_2-j_0-j_\infty)-\codim M_0^0(j_1,\ldots,j_k)=\\
&=2j_0+2j_\infty+\sum_{m=1}^{k} (2j_m-1),
\end{split}
\end{equation*}
which is the expected codimension.
\fi}
\end{proof}

\begin{rem}
If we would define the stratum $M(p)$ not inside $\C^{B_1\x}\x\C^{B_2\x}$, but inside $\C^{B_1}\x\C^{B_2}$, it would be wrong even in the classical case $B_i=\{0,1,2,\ldots,d_i\}$. Indeed, consider the stratum $M_0^0(\underbrace{1,\ldots,1}_{d_2})$. A general polynomial $f_2 \in \C^{B_2\x}$ has $d_2$ different roots, all of which are the roots of $f_1\equiv 0$, thus the subset $M_0^0(1,\ldots,1) \cap (\{0\} \x \C^{B_2\x})$ is open in $\{0\} \x \C^{B_2\x}$ and has codimension $(d_1+1)$. If $d_1 + 1 < d_2$, then the codimension of $M_0^0(1,\ldots,1) \cap (\{0\} \x \C^{B_2\x})$ is less than $d_2$, which is the expected codimension of $M_0^0(1,\ldots,1)$.
\end{rem}

The lemma survives as an estimate in the general case.
\begin{sledst}\label{corcodimMatmost}
In general, the codimension of the stratum $M(p)$ at every its irreducible component is at most the expected one.
\end{sledst}

Thus $M(p)$ can have irreducible components of larger dimension than expected, but can not have irreducible components of smaller dimension.

\begin{proof}
Equivalently, we can prove that the codimension of the stratum $M(p)$ at every its point is at most the expected one.

Consider the natural embedding $\C^{B_1\x}\times\C^{B_2\x}\to \C^{\conv B_1\x}\times\C^{\conv B_2\x}$, where $\conv B_i=\{ \min B_i, \ldots, \max B_i \}$ is the convex hull of $B_i$. The stratum $M(p)$ in $\C^{B_1\x}\times\C^{B_2\x}$ is the preimage of the corresponding stratum in $\C^{\conv B_1\x}\times\C^{\conv B_2\x}$. Thus at every point of the stratum $M(p)$ in $\C^{B_1\x}\times\C^{B_2\x}$ its codimension is not greater than the codimension of the corresponding stratum in $\C^{\conv B_1\x}\times\C^{\conv B_2\x}$ at the image point, which is equal to the expected one.
\end{proof}

The inverse estimate is in general not true, it requires some conditions on $B_1$ and $B_2$.

\begin{exa}
For example, if one can shift $B_1$ and $B_2$ to the same proper sublattice $k\Z\subset\Z$ with $k\geqslant 2$, then $\widehat M(1)=\widehat M(1,1)=\ldots=\widehat M(1,\ldots,1)$, where the last $\widehat M(1,\ldots,1)$ has $k$ ones. But the expected codimensions of $\widehat M(1)$ and $\widehat M(1,\ldots,1)$ are 1 and $k$ respectively, thus some of the strata $M(p)$ do not have the expected codimension (usually $M(1,\ldots,1)$ has codimension 1 instead of $k$). % TODO: уточнить
\end{exa}

% We will prove the inverse estimate in low codimensions, but it is better formulated in slightly different context. In this context, the inverse estimate is given for the codimensions up to 3, while for the strata $M(p)$ it is given only for the codimensions up to 2. The context is given in Subsection \ref{subsecintstr}, and the result is stated in Theorem \ref{thcodim}.

% But first let us study the strata $M(p)$ in more detail. % TODO: reformulate

\subsection{Proof of the part (ii) of the main theorem}\label{subsecclosstrat}
Using the methods of Lemma \ref{lemcodimclassic}, one can prove the similar facts about the closed strata.

\begin{lemma}
In the classical case $B_i=\{0,1,2,\ldots,d_i\}$, the set $\widehat M(p)$ is closed and its codimension of $\widehat M(p)$ is equal to the expected codimension (of $M(p)$).
\end{lemma}

{\iflong
\begin{proof}
Consider a subset $\Pr \widehat M(p) \subset \Pr(\C^{B_1}) \times \Pr(\C^{B_2})$ which is a projectivization of $\widehat M(p) \subset \C^{B_1\x}\x\C^{B_2\x}$. One can similarly to the proof of Lemma \ref{lemcodimclassic} construct a closed subset $\widehat P \subset (\CP^1)^{e_1} \times (\CP^1)^{e_2} \times (\CP^1)^k$ and a surjective map $\widehat P \to \Pr \widehat M(p)$ with finite fibers, thus $\Pr \widehat M(p)$ is closed and its codimension is equal to the expected one (at every point, and thus at every irreducible component). Moreover, the same holds for $\widehat M(p)$.
\end{proof}
\fi}

\begin{sledst}\label{corcodimwidehatM}
In general, the closed stratum $\widehat M(p)$ is indeed closed and its codimension (at every its irreducible component) is at most the expected one.
\end{sledst}

{\iflong
It can be deduced from the lemma above similarly to Corollary \ref{corcodimMatmost}.
\fi}

\begin{lemma}\label{lemclosm}
If $|B_1|>2$ and $B_2\neq\{j,j+1\}$, then $\{0\} \x \C^{B_2}$ lies in the closure of $\widehat M(1,1)$.
\end{lemma}

Notice that both conditions in the lemma are necessary. If $|B_1|=2$, then both $\{0\} \x \C^{B_2}$ and $\widehat M(1,1)$ have codimension 2 (and $\widehat M(1,1)$ may be empty). If $B_2=\{i,i+1\}$, then $\widehat M(1,1)$ is empty.

\begin{proof}
To show that $(0,g) \in \{0\} \x \C^{B_2}$ lies in this closure, it is enough to find $f \in \C^{B_1\x}$ such that $f$ and $g$ has two common roots, thus we would have that $(f, g) \in \widehat M(1,1)$ and $(0,g)=\lim_{t\to 0}(tf, g)$.

Suppose that $g$ has two distinct non-zero roots $x_1$ and $x_2$. The system
$$f(x_1)=\sum_{i \in B_1} f_i x_1^i=0 \ \text{ and } \ f(x_2)=\sum_{i \in B_1} f_i x_2^i=0$$
is a system of two homogeneous equations on $|B_2|>2$ variables, thus it has a non-zero solution $f$. This solution is a polynomial $f \in \C^{B_1\x}$ with two common roots with $g$ that is necessary to the proof.

If $B_2=\{j,j+1\}$ for some $i \in \Z$, then any $g \in \C^{B_2\x}$ has only one non-zero root. Otherwise a general polynomial $g \in \C^{B_2\x}$ has at least two different roots, thus $(0,g)$ lies in the closure of $\widehat M(1,1)$. As it is true for a general $g$, the subset $\{0\} \x \C^{B_2}$ also lies there. % TODO: details
\end{proof}

\begin{proof}[Proof of the main theorem \ref{theormain}, part (ii)]
We would like to have not only the condition (6) of the main theorem \ref{theormain}, but also the condition (1). Thus let us take the largest $k$ such that $B_1=k \cdot B_1'+m_1$ and $B_2=k \cdot B_2'+m_2$ for some $B_1'$ and $B_2' \subset \Z$: now we can not shift $B_1'$ and $B_2'$ to the same proper sublattice $l\Z\subset\Z$ with $l\geqslant 2$. As we know from Example \ref{exaintro}, part 1), $R_B \cong R_{B'}$.

By Lemma \ref{corcodimwidehatM} the subset $\widehat M(1,1)$ is closed in $\C^{B_1'\x}\x\C^{B_2'\x}$ and has codimension at most two (at every its irreducible component). But by Example \ref{exawidehat0011} $\widehat M(1,1)=\sing R_{B'} \cap (\C^{B_1'\x}\x\C^{B_2'\x})$. As $B_1'$ and $B_2'$ can not be shifted to the same proper sublattice $k\Z\subset\Z$ with $k\geqslant 2$, the subset $R_{B'}$ has codimension at least 1, thus the subset $\widehat M(1,1)=\sing R_{B'} \cap (\C^{B_1'\x}\x\C^{B_2'\x})$ has codimension at least 2. Thus $\widehat M(1,1)$ has codimension exactly 2.

Let us now notice that by Corollary \ref{corsing}
$$\sing R_{B'}=\begin{cases}
\widehat M(1,1) \cup \{(0,0)\} \quad \text{ if } \  B_1'=\{i,i+1\} \ \text { and } \ B_2'=\{j,j+1\} \\
\widehat M(1,1) \cup (\{0\} \x \C^{B_2'}) \quad \text{ if } \ B_1'=\{i,i+1\} \ \text { and } \ B_2'\neq\{j,j+1\} \\ 
\widehat M(1,1) \cup (\C^{B_1'} \x \{0\}) \quad \text{ if } \ B_1'\neq\{i,i+1\} \ \text { and } \ B_2'=\{j,j+1\} \\  
\widehat M(1,1) \cup (\{0\} \x \C^{B_2'}) \cup (\C^{B_1'} \x \{0\}) \quad \text{ otherwise } 
\end{cases}$$

If $|B_1'|=2$, then the subset $\{0\} \x \C^{B_2}$ has codimension 2, which is ok. Otherwise, we have $|B_1'|>2$ and $B_2' \neq \{j,j+1\}$, thus Lemma \ref{lemclosm} says that $\{0\} \x \C^{B_2'}$ lies in the closure of $\widehat M(1,1)$, so it does not give a component of $\sing R_{B'}$ of larger codimension.

The same can be proven about $\C^{B_1'} \x \{0\}$.

As for the subset $\{(0,0)\}$, it appears by itself only in the case $B_1'=\{i,i+1\}$ and $B_2'=\{j,j+1\}$. Then $B_1=k \cdot B_1'+m_1=\{i',i'+k\}$ and $B_2=k \cdot B_2'+m_2=\{j',j'+k'\}$, which is excluded by the condition (6) of the main theorem \ref{theormain}.
\end{proof}

\section{A stratification of the solution space}\label{secstratrootspace}
In this section, we will prove Theorem \ref{theorcodimNpre} modulo Lemma \ref{lemstratumranks11}, Theorem \ref{theorstratumranks21for111} and Lemma \ref{lem3minor2} from the subsequent sections (the proofs of Lemmas \ref{lemcork2simple} and \ref{lemcork2} refer to these facts), thus finishing the proof of the main theorem \ref{theormain}.

\subsection{The relation to the stratication of the sparse resultant}

Suppose that $j_0=j_\infty=0$ and study the filtration subsets $N\left(\,^{j^1_1,\ldots,j^1_k}_{j^2_1,\ldots,j^2_k}\right) \subset \C^{B_1\x} \times \C^{B_2\x}$ (in particular, $N(j_1,\ldots,j_k)=N\left(\,^{j_1,\ldots,j_k}_{j_1,\ldots,j_k}\right)$; see Definition \ref{defNN}). Fix $B_1$ and $B_2 \subset \Z$ and a subset $N\left(\,^{j^1_1,\ldots,j^1_k}_{j^2_1,\ldots,j^2_k}\right)$.

Consider the solution space
$$Z_k=\left\{ (x_1, \ldots, x_k) \subset (\C^\times)^k \,|\, x_i \neq x_j \; \text{for} \; i \neq j \right\}$$
of tuples of different non-zero numbers.

We would like to define a stratification on it that is related with the structure of $N\left(\,^{j^1_1,\ldots,j^1_k}_{j^2_1,\ldots,j^2_k}\right)$.

\begin{defin}\label{defmatrsimple}
First let us for each tuple $x=(x_1,\ldots,x_k)$ define a pair of polynomials
$$h^1_x(t)=\prod_{m=1}^k (t-x_m)^{j_m^1} \quad \text{ and } \quad h^2_x(t)=\prod_{m=1}^k (t-x_m)^{j_m^2}.$$

Then consider a pair of maps
$$\psi^1_x: \C^{B_1} \to \C[t]/(h^1_x(t)) \quad \text{ and } \quad \psi^2_x: \C^{B_2} \to \C[t]/(h^2_x(t))$$
given by
$$f(t) \mapsto f(t) \mod h^1_x(t) \quad \text{ and } \quad g(t) \mapsto g(t) \mod h^2_x(t).$$

(c.f. \cite{E13}*{Def. 3.21}).

We can now define the strata
$$S^1_{n_1} = \left\{ (x_1, \ldots, x_k) \subset Z_k \,|\, \codim \Im \psi^1_x = n_1 \right\},$$
$$S^2_{n_2} = \left\{ (x_1, \ldots, x_k) \subset Z_k \,|\, \codim \Im \psi^2_x = n_2 \right\},$$
$$S_{n_1, n_2} = S^1_{n_1} \cap S^2_{n_2} = \left\{ (x_1, \ldots, x_k) \subset Z_k \,|\, \codim \Im \psi^1_x = n_1 \text{ and } \codim \Im \psi^2_x = n_2 \right\}.$$

To sum it up, we fix $B_1$ and $B_2 \subset \Z$, fix a subset $N\left(\,^{j^1_1,\ldots,j^1_k}_{j^2_1,\ldots,j^2_k}\right)$ and define a stratification
$$Z_k = \bigsqcup_{n_1, n_2} S_{n_1, n_2} \subset (\C^\times)^k.$$
\end{defin}

\begin{rem}\label{remM}
The polynomial $h^i_x(t)$ vanishes on $x_m$: $h^i_x(x_m)=0$; moreover, a number of its derivatives also vanish: $(h^i_x)^{(d^i_m)}(x_m)=0$ for $d^i_m<j^i_m$; as a consequence, the space $\C[t]/(h^i_x(t))$ has a system of coordinate functions sending the class $f(t) \mod h^i_x(t)$ to
$$f(x_1), f'(x_1), \ldots, f^{(j^i_1-1)}(x_1), \quad f(x_2), f'(x_2), \ldots, f^{(j^i_2-1)}(x_2), \quad \ldots, \quad f(x_k), f'(x_k), \ldots, f^{(j^i_k-1)}(x_k).$$
These functions form a basis of the dual space. Using this basis, we can write down the matrix $M_i$ of the map $\psi^i_x$.

For example, suppose that $N\left(\,^{j^1_1,\ldots,j^1_k}_{j^2_1,\ldots,j^2_k}\right)=N(1,1,1)$. Then the matrix of the map $\psi^i_x$ has the form
\begin{equation}\label{eqM}
M_i=M_i(x,y,z) = \begin{pmatrix}
x^{b^i_1} & x^{b^i_2} & x^{b^i_3} & \ldots & x^{b^i_{s_i}} \\
y^{b^i_1} & y^{b^i_2} & y^{b^i_3} & \ldots & y^{b^i_{s_i}} \\
z^{b^i_1} & z^{b^i_2} & z^{b^i_3} & \ldots & z^{b^i_{s_i}}
\end{pmatrix}
\end{equation}
where the columns are numbered by the elements of $B_i$ and $x_1,x_2,x_3$ are replaced by $x,y,z$.

Moreover, we have
$$\codim \Im \psi^i_x=\cork M_i,$$
where $M_i$ is some matrix depending on $x_1,\ldots,x_k$. 
\end{rem}

% The matrices $M_i = M(B_i; \ j^i_1, \ldots, j^i_k; \ x_1, \ldots, x_k)$ have $|J_i|= \sum_{m=1}^{k} j^i_m$ rows and $|B_i|$ columns respectively. Their coranks are the numbers $\cork M_i=|J_i|-\rk M_i$, that is the numbers of their rows minus their rank.

\begin{rem}
Now consider the subsets
$$\widehat S^i_{n_i} = \bigcup_{m_i \geqslant n_i} S^i_{m_i} = \left\{ (x_1, \ldots, x_k) \subset Z_k \,|\, \codim \Im \psi^i_x \geqslant n_i \right\}.$$
By the previous remark,
$$\widehat S^i_{n_i}= \left\{ (x_1, \ldots, x_k) \subset Z_k \,|\, \cork M_i \geqslant n_i \right\},$$ 
thus they are closed. Consequently, the strata $S^i_{n_i}$ and $S_{n_1, n_2}=S^1_{n_1} \cap S^2_{n_2}$ are locally closed.
\end{rem}

Thus theorem \ref{theorcodimNpre} reduces to estimating the codimension of $S_{n_1,n_2}$.

Our interest in these strata is due to the following lemma:
\begin{lemma}\label{lemcodimS}
If for each $S_{n_1, n_2}$ the codimension (of every its irreducible component) is at least $n_1 + n_2$, then the codimension of (every irreducible component of) $N\left(\,^{j^1_1,\ldots,j^1_k}_{j^2_1,\ldots,j^2_k}\right)$ is at least the expected codimension $\sum_{m=1}^k(j^1_m+j^2_m-1)$.
\end{lemma}

Typically, we would have a simple set of conditions implying that $\codim S^1_{n_1} \geqslant n_1$ and $\codim S^2_{n_2} \geqslant n_2$, and the difficult part would be to estimate the codimension of their intersection.

\begin{proof}
Consider the subset $\widetilde N=\widetilde N\left(\,^{j^1_1,\ldots,j^1_k}_{j^2_1,\ldots,j^2_k}\right)$ consisting of $(f_1;\ f_2; x_1, \ldots, x_k) \in \C^{B_1\x} \times \C^{B_2\x} \times Z_k$ such that for $m=1,\ldots,k$ holds $\ord_{x_m}(f_1) \geqslant j^1_m$ and $\ord_{x_m}(f_2) \geqslant j^2_m$.

The image of $\widetilde N$ under the projection $p: \C^{B_1\x} \times \C^{B_2\x} \times Z_k \to \C^{B_1\x} \times \C^{B_2\x}$ is exactly $N=N\left(\,^{j^1_1,\ldots,j^1_k}_{j^2_1,\ldots,j^2_k}\right)$, and the preimage of any point of $N$ consists of a finite number of points. We have that $\dim\ (\C^{B_1\x} \times \C^{B_2\x} \times Z_k) = \dim\ (\C^{B_1\x} \times \C^{B_2\x}) + k$. Thus if we prove that the codimension of (every irreducible component of) $\widetilde N$ is at least $\sum_{m=1}^k(j^1_m+j^2_m)$, we will prove that the codimension of (every irreducible component of) $N$ is at least the expected codimension $\sum_{m=1}^k(j^1_m+j^2_m-1)$.

Now consider the projection $q: \C^{B_1\x} \times \C^{B_2\x} \times Z_k \to Z_k$. The fiber over the point $(x_1, \ldots, x_k) \in Z_k$ is a subset in $\C^{B_1\x} \times \C^{B_2\x}$ of all pairs of non-zero polynomials $(f_1, f_2)$ such that at the point $x_m$ the polynomial $f_1$ has a root of multiplicity at least $j^1_m$ and the polynomial $f_2$ has a root of multiplicity at least $j^2_m$. It is the direct product of the spaces
$$\ker \psi^1_x \setminus \{0\} = \left\{ f_1 \in \C^{B_1\x} \,|\, \ord_{x_m} f_1 \geqslant j^1_m \right\} \text{ and } \ker \psi^2_x \setminus \{0\}  = \left\{ f_2 \in \C^{B_2\x} \,|\, \ord_{x_m} f_2 \geqslant j^2_m \right\}.$$

Let us notice that
$$\codim (\ker \psi^i_x \setminus \{0\}) \geqslant \codim \ker \psi^i_x=\dim \Im \psi^i_x=\sum_{m=1}^{k} j^i_m-\codim \Im \psi^i_x=\sum_{m=1}^{k} j^i_m-n_i,$$
thus for the point $(x_1, \ldots, x_k) \in S_{n_1, n_2} \subset Z_k$ the codimension of the fiber
$$q^{-1}(x_1, \ldots, x_k)=(\ker \psi^1_x \setminus \{0\}) \times (\ker \psi^2_x \setminus \{0\}) \subset \C^{B_1\x} \times \C^{B_2\x}$$
is at least $\sum_{m=1}^{k} (j^1_m+j^2_m) - n_1 - n_2$.

If the codimension of (every irreducible component of) $S_{n_1,n_2}$ is at least $n_1+n_2$, then the codimension of (every irreducible component of) $\widetilde N \cap q^{-1}(S_{n_1,n_2})$ is at least $\sum_{m=1}^{k} (j^1_m+j^2_m)$. Consequently, the codimension of (every irreducible component of) $\widetilde N$ is at least $\sum_{m=1}^{k} (j^1_m+j^2_m)$, and the codimension of (every irreducible component of) $N$  is at least $\sum_{m=1}^{k} (j^1_m+j^2_m-1)$, as necessary.
\end{proof}

\subsection{The stratification for one polynomial}

Let us fix $i=1$ or $2$ and first consider only one of the stratifications $S^i_{n_i}$. We will take a subset $N\left(\,^{j^1_1,\ldots,j^1_k}_{j^2_1,\ldots,j^2_k}\right) \subset \C^{B_1\x} \times \C^{B_2\x}$ and study the corresponding stratifications $Z_k = \bigsqcup_{n_i} S^i_{n_i}$ for small $k$.

We wonder how the structure of $Z_k = \bigsqcup_{n_i} S^i_{n_i}$ depends on $N\left(\,^{j^1_1,\ldots,j^1_k}_{j^2_1,\ldots,j^2_k}\right)$ and $B_i$. It depends only on the coefficients $j^i_1,\ldots,j^i_k$ for the fixed $i$, not on the whole $N\left(\,^{j^1_1,\ldots,j^1_k}_{j^2_1,\ldots,j^2_k}\right)$, thus we will simply write "for $(j^i_1,\ldots,j^i_k)$" instead of "for $N\left(\,^{j^1_1,\ldots,j^1_k}_{j^2_1,\ldots,j^2_k}\right)$".

To simplify the notations, we need a following definition:

\begin{defin}
For $B_i=\{ b^i_1, \ldots, b^i_k \}$ let us define
$$\phi(B_i) := (b^i_2-b^i_1,\ b^i_3-b^i_2,\ \ldots,\ b^i_k-b^i_{k-1})$$
to be the greatest common divisor of the differences.
\end{defin}

Thus the condition (1) of the main theorem \ref{theormain} can be written as $(\phi(B_1), \phi(B_2))\geqslant 2$, where $(\phi(B_1), \phi(B_2))$ is the greatest common divisor of $\phi(B_1)$ and $\phi(B_2)$.

First we will study the case of $(1,1,1)$. Here we use the matrices $M_i$ from Formula (\ref{eqM}) of Remark \ref{remM}.

\begin{lemma}\label{lemcork1simple}
The subset $S^i_2$ for $(1,1,1)$ is non-empty if and only if $\phi(B_i)\geqslant 3$.

Namely, all such triples $(x,y,z)$ that $\cork M_i(x,y,z)=2$ have the form $(c,ct,cu)$, where $t$ and $u$ are two $\phi(B_i)$-th roots of unity such that $t \neq u$, $t \neq 1$ and $u \neq 1$.

In particular, if $S^i_2$ is non-empty, it has codimension 2.
\end{lemma}

\begin{proof}
If $c$ is a non-zero number, then the matrices $M_i(x,y,z)$ and $M_i(cx,cy,cz)$ differ by the multiplication of their $j$-th columns on non-zero numbers $c^{b_j^i}$, thus their ranks are the same and $(x,y,z) \in S^i_{n_i}$ if and and only if $(cx,cy,cz) \in S^i_{n_i}$.

Thus we can replace $M_i(x,y,z)$ with $M_i(1,t,u)$, where $t=y/x$ and $u=z/x$. The matrix
$$M_i(1,t,u)=\begin{pmatrix}
1 & 1 & 1 & ... & 1 \\ 
t^{b^i_1} & t^{b^i_2} & t^{b^i_3} & ... & t^{b^i_{s_i}} \\
u^{b^i_1} & u^{b^i_2} & u^{b^i_3} & ... & u^{b^i_{s_i}}
\end{pmatrix}$$
has corank 2 if and only if all its $2 \times 2$-minors are degenerate, thus
$$t^{b_1^i}=t^{b_2^i}=t^{b_3^i}=\ldots=t^{b_{s_i}^i} \quad \text{ and } \quad  u^{b_1^i}=u^{b_2^i}=u^{b_3^i}=\ldots=u^{b_{s_i}^i}.$$
But $t \neq 0$, thus $t$ is a $k_1$-th root of unity for some $k_1$ and all $(b^i_q-b^i_p)$ should be divisible by $k_1$. Similarly, $u$ is a $k_2$-th root of unity for some $k_2$ and all $(b^i_q-b^i_p)$ should be divisible by $k_2$. Now $t$ and $u$ are both $k$-th roots of unity for the greatest common divisor $k=(k_1, k_2)$, and all $(b^i_q-b^i_p)$ should be divisible by $k$. Moreover, $t \neq u$, $t \neq 1$ and $u \neq 1$, thus there are at least three different $k$-th roots of unity, thus $k \geqslant 3$. To sum it up, $\phi(B_i) \geqslant 3$.

Vice versa, if $\phi(B) = k \geqslant 3$, then we can choose three different $k$-th roots of unity, namely 1, $t$ and $u$, and the matrix $M_i(1,t,u)$ will be of the corank 2.
\end{proof}

\begin{sledst}\label{corcodim111}
Consider a decomposition $Z_3=S^i_0 \cup S^i_1 \cup S^i_2$ for $(1,1,1)$:
\begin{itemize}
    \item If $|B_i| \geqslant 3$, then holds $Z_3=\begin{cases}
    S^i_0 \cup S^i_1 \cup S^i_2 \text{ with } \codim S^i_1 \geqslant 1 \\
    \quad \quad \quad \quad \quad \text { and } S^i_2 \text{ non-empty of codim 2 if } \phi(B_i) \geqslant 3 \\
    S^i_0 \cup S^i_1 \text{ with } \codim S^i_1 \geqslant 1 \text{ otherwise }
    \end{cases}$
    \item If $|B_i|=2$, then holds $Z_3=\begin{cases}
    S^i_1 \cup S^i_2 \text{ with } S^i_2 \text{ non-empty of codim 2 if } \phi(B_i) \geqslant 3 \\
    S^i_1 \text { otherwise }
    \end{cases}$
\end{itemize}
\end{sledst}

\begin{proof}
The subset $S^i_0$ consists of $(x,y,z)$ such that the matrix $M_i(x,y,z)$ from Formula (\ref{eqM}) of Remark \ref{remM} has a non-degenerate $3 \times 3$-minor.

Suppose that $|B_i|\geqslant 3$: then the subset $S^i_0$ is open. Its complement $S^i_1 \cup S^i_2$ is a subset of codimension at least 1, but by Lemma \ref{lemcork1simple} the subset $S^i_2$ is of codimension 2 (or empty), thus $S^i_1$ is also of codimension at least 1. Moreover, by Lemma \ref{lemcork1simple}, the subset $S^i_2$ is non-empty if and only if $\phi(B_i)\geqslant 3$.

Now suppose that $|B_i|=2$: then there is no $3 \times 3$-minors in $M_i$, thus the subset $S^i_0$ is empty. Thus $Z_3=S^i_1 \cup S^i_2$, but $S^i_2$ is a subset of codimension 2 (or empty), thus now $S^i_1$ is open.
\end{proof}

Now let us consider the other filtration subsets $N(j_1,\ldots,j_k)$ of expected codimension 1, 2 or 3. The other cases are proven relatively straightforward compared to $N(1,1,1)$, which is why we devoted a separate statement to $N(1,1,1)$. We have a following lemma (cf. \cite{N19}*{Th. 1.1}): % TODO: check ref

\begin{lemma}\label{lemcork1}
i) For $(k)$ there are following equations:
\begin{itemize}
    \item If $|B_i| \geqslant k$, then holds $Z_1=S^i_0$.
    \item If $|B_i| \leqslant k-1$, then holds $Z_1=S^i_{k-|B_i|}$
\end{itemize}

ii) For $(k,1)$ there are following decompositions:
\begin{itemize}
    \item If $|B_i|\geqslant k+1$, then holds $Z_2=S^i_0 \cup S^i_1$ with $\codim S^i_1 \geqslant 1$
    \item If $|B_i|\leqslant k$, then holds $Z_2=S^i_1$
\end{itemize}

Moreover, for $(1,1)$ we have that $S^i_1$ is non-empty if and only if $\phi(B_i) \geqslant 2$.
\end{lemma}

In particular, this lemma covers the cases of $N(1)$, $N(2)$ and $N(1,1)$, but it also covers the case of $N\left(\,^{2,1}_{1,1}\right)$.

\subsection{The stratification for two polynomials}

Now let us fix a subset $N\left(\,^{j^1_1,\ldots,j^1_k}_{j^2_1,\ldots,j^2_k}\right) \subset \C^{B_1} \times \C^{B_2}$ and consider the stratification $Z_k = \bigsqcup_{n_1,n_2} S_{n_1,n_2}$ with $S_{n_1, n_2}=S^1_{n_1} \cap S^2_{n_2}$. We wonder how its structure depends on $N\left(\,^{j^1_1,\ldots,j^1_k}_{j^2_1,\ldots,j^2_k}\right)$ and $(B_1, B_2)$.

\begin{lemma}\label{lemcork2simple}
Suppose that $|B_1|\geqslant 3$ and $|B_2|\geqslant 3$. For $N(1,1,1)$ there is a decomposition
$$Z_3=\bigsqcup_{\genfrac{}{}{0pt}{2}{0\leqslant i\leqslant 2}{0\leqslant j\leqslant 2}} S_{ij}$$
such that
\begin{itemize}
    \item $S_{00}$ is open,
    \item $\codim S_{01}$ and $\codim S_{10} \geqslant 1$,
    \item $\codim S_{02}$ and $\codim S_{20} \geqslant 2$,
    \item if $B_1 \setminus \{\max(B_1)\}$ and $B_2 \setminus \{\max(B_2)\}$ can \textbf{not} be shifted to the same proper sublattice $k\Z\subset\Z$ with $k\geqslant 2$, and $B_1 \setminus \{\min(B_1)\}$ and $B_2 \setminus \{\min(B_2)\}$ can \textbf{not} be shifted to the same proper sublattice $k\Z\subset\Z$ with $k\geqslant 2$, then $\codim S_{11} \geqslant 2$,
    \item if $B_1$ can be split into $B'\sqcup B''$ so that $B'$, $B''$ and $B_2$ can be shifted to the same sublattice $k\Z\subset\Z$ with $k\geqslant 3$, but $B_1$ and $B_2$ can \textbf{not} be shifted to the same sublattice $k\Z$, then $\codim S_{12} = 2$ and it is non-empty, otherwise it is empty (and symmetrically with $S_{21}$), % TODO: упростить через \widehat S
    \item if $B_1$ and $B_2$ can be shifted to the same proper sublattice $k\Z\subset\Z$ with $k\geqslant 3$, then $\codim S_{22} = 2$ and it is non-empty, otherwise it is empty,
\end{itemize}
and all the strata can be empty if not mentioned otherwise.
\end{lemma}

To simplify the proof, we can use homogeneity considerations, that is replace
$$Z_3=\left\{ (x, y, z) \subset (\C^\times)^3 \,|\, x \neq y \neq z \neq x \right\}$$
with
$$\Pr Z_3=\left\{ (1:t:u) \subset \CP^2 \,|\, t,u \neq 0 \text{ and } t \neq u, t \neq 1, u \neq 1 \right\},$$
and replace the strata $S^i_{n_i}$ (and $S_{n_1,n_2}$) with their projectivizations $\Pr S^i_{n_i}$ (and $\Pr S_{n_1,n_2}$ respectively), having the same codimensions of components. % If the codimension (of every irreducible component) of $S^i_{n_i}$ (or $ S_{n_1, n_2}$) is equal/less than/greater than some number, then the same is true about the codimension (of every irreducible component) of $\Pr S^i_{n_i}$ (or $\Pr S_{n_1, n_2}$ respectively).

\begin{proof}
The cases of $S_{00}$, $S_{01}$, $S_{10}$, $S_{02}$ and $S_{20}$ follow from Corollary \ref{corcodim111}. 

\medskip
The case of $S_{11}=S_1^1 \cap S_1^2$ % for $N(1,1,1)$
follows from Lemma \ref{lemstratumranks11} of Section \ref{secstrats11forn111}. Indeed, by Lemma \ref{lemcork1simple} all the irreducible components of $\Pr\widehat S_1^1=\Pr S_1^1 \ \cup \ \Pr S_2^1$ and of $\Pr\widehat S_1^2$ have codimension at least 1. If $\Pr S_1^1$ has an irreducible component of codimension 1, then $\Pr\widehat S_{11} = \Pr\widehat S_1^1 \cap \Pr\widehat S_1^2 = \Pr S_{11} \cup \Pr S_{12} \cup \Pr S_{21} \cup \Pr S_{22}$ has an irreducible component $W \subset \Pr Z_3$ of codimension 1. The subset $\Pr\widehat S_{11} \subset \Pr Z_3$ is defined by $3 \times 3$-minors of $M_1$ and $M_2$, that is by the polynomials
$$\det\nolimits_{a,b,c}(t,u)=\begin{pmatrix}
1 & 1 & 1 \\ 
t^a & t^b & t^c\\
u^a & u^b & u^c
\end{pmatrix},$$
where we take the triples of elements $\{a, b, c\} \subset B_1$ or $\{a, b, c\} \subset B_2$ such that $a < b < c$. Let $S$ be a set of such polynomials.

We would like apply Lemma \ref{lemstratumranks11} to $S$, where $f_{a,b,c}(t,u)\equiv 0$ for each $\det_{a,b,c}+f_{a,b,c} \in S$. The polynomial $GCD(S)$ vanishes on $\Pr\widehat S_{11} \subset \Pr Z_3$, where $\Pr Z_3$ is a complement of $V(tu(t-1)(u-1)(t-u)) \subset \C^2$, and $\Pr\widehat S_{11}$ is non-empty, thus $GCD(S)$ does not divide a power of $tu(t-1)(u-1)(t-u)$.

Suppose that the first alternative of Lemma \ref{lemstratumranks11} holds: $k=GCD(\{b-a  \,|\, \det\nolimits_{a,b,c}+f_{a,b,c} \in S\}) > 1$. Then all the differences $(b^1_i-b^1_j)$ where $b^1_i \neq \max B_1$ and $b^1_j \neq \max B_1$ are divisible by $k$, and all the differences $(b^2_i-b^2_j)$ where $b^2_i \neq \max B_2$ and $b^2_j \neq \max B_2$ are divisible by $k$, thus $B_1 \setminus \{\max(B_1)\}$ and $B_2 \setminus \{\max(B_2)\}$ can be shifted to the same proper sublattice $k\Z\subset\Z$ with $k\geqslant 2$. % TODO:

Suppose that the second alternative of Lemma \ref{lemstratumranks11} holds: similarly, $B_1 \setminus \{\min(B_1)\}$ and $B_2 \setminus \{\min(B_2)\}$ can be shifted to the same proper sublattice $k\Z\subset\Z$ with $k\geqslant 2$.

\medskip
The case of $S_{12}=S_1^1 \cap S_2^2$ (and similarly $S_{21}=S_2^1 \cap S_1^2$) % for $N(1,1,1)$
is the most difficult part of the lemma. It follows from Theorem \ref{theorstratumranks21for111} of Section \ref{secstrats12forn111}. Indeed, by Lemma \ref{lemcork1simple} the codimension of $S^2_2$ is 2 (maybe empty), thus the codimension of $S_{12}$ is at least 2, while it is at most 2, because it is equal to the codimension of $\Pr S^2_2$. Thus we should check the conditions when $S_{12}$ is non-empty. % TODO: delete most difficult here and everywhere

If $S_{12}$ is non-empty, then we can take $(1:t:u) \in \Pr S_{12}$ and apply Lemma \ref{lemcork1simple}. As $(1:t:u) \in \Pr S_2^2$, then $t$ and $u$ are $n$-th roots of unity, where $n=\phi(B_2)$. Moreover, $(1:t:u) \in \Pr \widehat S_1^1$, thus we can apply Theorem \ref{theorstratumranks21for111} to $x=t$, $y=u$ and $B=B_1$ and get that $B_1$ can be split into $B'\sqcup B''$ so that $\phi(B')$, $\phi(B'')$ and $n=\phi(B_2)$ have a common divisor $k\geqslant 3$. It means exactly that $B'$, $B''$ and $B_2$ can be shifted to the same sublattice $k\Z\subset\Z$ with $k\geqslant 3$, as necessary. Moreover, suppose that $B_1$ and $B_2$ can be shifted to the same proper sublattice $k\Z\subset\Z$. Then $(1:t:u) \in \Pr S_{22}$, thus $(1:t:u) \not\in \Pr S_{12}$, a contradiction.

Vice versa, if $B_1$ can be split into $B'\sqcup B''$ so that $B'$, $B''$ and $B_2$ can be shifted to the same sublattice $k\Z\subset\Z$ with $k\geqslant 3$, then we can choose three different $k$-th roots of unity, namely 1, $t$ and $u$. As $k$ divides $\phi(B_1)$, we have that $(1:t:u) \in \Pr S^2_2$ by Lemma \ref{lemcork1simple}. The subset $\Pr \widehat S_1^1$ is defined by $3 \times 3$-minors of $M_1$, that is by the polynomials
$$\det\nolimits_{a,b,c}(t,u)=\begin{pmatrix}
1 & 1 & 1 \\ 
t^a & t^b & t^c\\
u^a & u^b & u^c
\end{pmatrix},$$
where we take the triples of elements $\{a, b, c\} \subset B_1$. As $B_1 = B'\sqcup B''$, there are two of $a$, $b$ and $c$ from the same subset. Without loss of generality we can suppose that $a$ and $b \in B'$. As $k$ divides $\phi(B')$, we have that $t^a=t^b$ and $u^a=u^b$, thus the $3 \times 3$-minor $\det\nolimits_{a,b,c}$ has two equal columns and is zero. To sum it up, $(1:t:u) \in \Pr \widehat S^1_1$ and $(1:t:u)\in \Pr \widehat S_{12}$, thus $\widehat S_{12}=S_{12} \cup S_{22}$ is non-empty. 

Moreover, we have the condition that $B_1$ and $B_2$ can not be shifted to the same proper sublattice $k\Z\subset\Z$. Let $l$ be the largest number such that $B_1$ and $B_2$ can be shifted to the same proper sublattice $l\Z\subset\Z$. As $B_1=B'\sqcup B''$, we have that $k$ is divisible by $l$. Moreover, $k > l$, thus we can choose $t$ and $u$ in such a way that they are not both $l$-th roots of unity. It means that $(1:t:u)$ lies not in $S_{22}$, but in $S_{12}$, and thus $S_{12}$ is non-empty.

\medskip
The case of $S_{22}=S_2^1 \cap S_2^2$ % for $N(1,1,1)$
is proven as follows. By Lemma \ref{lemcork1simple}, $\Pr S_2^i$ consists of $(1:t:u)$, where $t$ and $u$ are two $\phi(B_i)$-th roots of unity such that $t \neq u$, $t \neq 1$ and $u \neq 1$. Thus if $(1:t:u) \in \Pr S_{22}=\Pr S_2^1 \cap \Pr S_2^2$, then $t$ and $u$ are $\phi(B_1)$-th roots of unity and $\phi(B_2)$-th roots of unity, thus they are $k=(\phi(B_1),\phi(B_2))$-th roots of unity. Let us recall that $t \neq 1$, $u \neq 1$ and $t \neq u$, thus $k \geqslant 3$. In other words, $B_1$ and $B_2$ can be shifted to the same proper sublattice $k\Z\subset\Z$ with $k\geqslant 3$. Both $\Pr S_2^1$ and $\Pr S_2^2$ has codimension 2 in the space $\Pr Z_3$ of dimension 2, thus $\Pr S_{22}$ is also of codimension 2.

Vice versa, if $B_1$ and $B_2$ can be shifted to the proper sublattice $k\Z\subset\Z$ with $k\geqslant 3$, we can take a number $t$ which is not equal to 1 and is a $k$-th root of unity, then $(1:t) \in \Pr S_{11}$ and thus $S_{11}$ is non-empty.
\end{proof}

Now we can prove the most difficult part of Theorem \ref{theorcodimNpre}.

\begin{sledst}\label{codimNsimple}
Under the conditions of the part (i) of the main theorem \ref{theormain}, the codimension of (every irreducible component of) $N(1,1,1)$ is at least the expected codimension 3.
\end{sledst}

\begin{proof}
By the condition (5) of the main theorem, we have that either $|B_1|=2$ and $\max B_2-\min B_2\leqslant 2$ (in which case $N(1,1,1)$ is empty), or $|B_1|>2$, and the same with $B_2$. % TODO: details?

Thus we can use Lemma \ref{lemcork2simple} and prove that $\codim S_{n_1,n_2}\geqslant n_1+n_2$. Indeed, the condition for $S_{11}$ follows from the conditions (3) and (4) of the main theorem, while the condition for $S_{12}$ and $S_{21}$ follows from the condition (2) of the main theorem.

Now by Lemma \ref{lemcodimS} we have $\codim N(1,1,1) \geqslant 3$.
\end{proof}

Now let us consider the other filtration subsets $N(j_1,\ldots,j_k)$ of expected codimension 1, 2 or 3.

\begin{lemma}\label{lemcork2}
There are following decompositions:
\begin{itemize}
    \item For $N(k)$ holds $Z_1=S_{00}$.
    \item For $N(1,1)$ holds $Z_2=S_{00} \cup S_{01} \cup S_{10} \cup S_{11}$, where
    \begin{enumerate}
        \item $\codim S_{01} = \codim S_{10} = 1$ (and they may be empty),
        \item if $B_1$ and $B_2$ can be shifted to the same proper sublattice $k\Z\subset\Z$ with $k\geqslant 2$, then $\codim S_{11} = 1$ and it is non-empty, otherwise it is empty.
    \end{enumerate}
    \item If $|B_1|>2$, then for $N\left(\,^{2,1}_{1,1}\right)$ holds $Z_3=S_{00} \cup S_{01} \cup S_{10} \cup S_{11}$, where
    \begin{enumerate}
        \item $\codim S_{01}$ and $\codim S_{10} \geqslant 1$ (and they may be empty),
        \item if $B_1$ and $B_2$ can be shifted to the same proper sublattice $k\Z\subset\Z$ with $k\geqslant 2$, then $\codim S_{11} = 1$ and it is non-empty, otherwise it is empty.
    \end{enumerate} % TODO: check
    \item If $|B_1|=2$, then for $N\left(\,^{2,1}_{1,1}\right)$ holds $Z_3= S_{10} \cup S_{11}$, where $\codim S_{11}=1$.  % TODO: check
\end{itemize}
\end{lemma}

Here we list those symmetric filtration subsets $N(j_1,\ldots,j_k)$ whose expected codimension is 1, 2 or 3, and also $N\left(\,^{2,1}_{1,1}\right)$. We write $S_{ij}$ instead of $S_{i,j}$ for shortness.

\begin{proof}
Most of the cases trivially follow from Lemmas \ref{lemcork1simple} and \ref{lemcork1}.

\medskip
The case of $S_{11}=S_2^1 \cap S_2^1$ for $N(1,1)$ can be proven similarly to the case of $S_{22}$ for $N(1,1,1)$.

\medskip
The case of $S_{11}=S_1^1 \cap S_1^2$ for $N\left(\,^{2,1}_{1,1}\right)$ follows from Lemma \ref{lem3minor2} similarly to how the case of $S_{11}$ for $N(1,1,1)$ follows from Lemma \ref{lem3minor}.
\end{proof}

Now we can formulate Theorem \ref{theorcodimNpre} in full detail:

\begin{theor}\label{theorcodimN}
The following strata have codimensions (of every their irreducible components) more then or equal to their expected codimensions unless the following conditions of the main theorem \ref{theormain} hold:
\begin{itemize}
    \item $\codim N(1) \geqslant 1$ always.
    \item $\codim N(1,1) \geqslant 2$ unless (1).
    \item $\codim N^0_1 \geqslant 2$ always.
    \item $\codim N^1_0 \geqslant 2$ always.
    \item $\codim N(2) \geqslant 3$ always.
    \item $\codim N(1,1,1) \geqslant 3$ unless (1)-(5).
    \item $\codim N_0^1(1) \geqslant 3$ always.
    \item $\codim N_1^0(1) \geqslant 3$ always.
    \item $\codim N\left(\,^{2,1}_{1,1}\right) \geqslant 3$ unless (1).
    \item $\codim N\left(\,^{1,1}_{2,1}\right) \geqslant 3$ unless (1).
\end{itemize}
\end{theor}

It is deduced from Lemmas \ref{lemcork2simple} and \ref{lemcork2} directly case-by-case.

\section{The stratum S(1,1) of the solution space for the filtration subset N(1,1,1)}\label{secstrats11forn111}

In this section we will study the stratum $S_{1,1}$ for $N(1,1,1)$, although without refering to it as such. This section is independent from the other parts of the text.

We will work inside the space
$$\Pr Z_3=\left\{ (t, u) \subset (\C^\x)^2 \,|\, t \neq 1, \ u \neq 1, \ t \neq u \right\}.$$

\begin{defin}
Consider a triple of integers $a<b<c$. It defines a (Laurent) polynomial
$$\det\nolimits_{a,b,c}(t,u)=\det \begin{pmatrix}
1 & 1 & 1 \\ 
t^a & t^b & t^c\\
u^a & u^b & u^c
\end{pmatrix}$$
for $a < b < c$. Its support is a hexagon $N_{det_{a,b,c}}$ with vertices $(a,b)$, $(a,c)$, $(b,c)$, $(c,b)$, $(c,a)$ and $(b,a)$, see Figure \ref{fignewtdet}.

Now take a polynomial $f_{a,b,c}(t,u)$ such that its support lies inside $N_{det_{a,b,c}}$ and such that it is divisible by $det_{0,1,2}$. We will call such a polynomial \textit{good} polynomial.
\end{defin}

\begin{lemma}\label{lemstratumranks11}
Let $S$ be a (not necessary finite) set of good polynomials $\det_{a,b,c}+f_{a,b,c}$. Suppose that their greatest common divisor $GCD(S)$ does not divide a power of $tu(t-1)(u-1)(t-u)$. Then at least one of the numbers
$$GCD(\{b-a  \,|\, \det\nolimits_{a,b,c}+f_{a,b,c} \in S\}) \quad \text{ and } \quad GCD(\{c-b  \,|\, \det\nolimits_{a,b,c}+f_{a,b,c} \in S\})$$
is greater than 1.
\end{lemma}

\begin{proof}
Let us decompose $GCD(S)$ as $g(t,u) \cdot (t-1)^\alpha \cdot (u-1)^\beta \cdot (t-u)^\gamma$, where $g(t,u)$ is not divisible by $(t-1)$, $(u-1)$ or $(t-u)$. Here we will consider divisibility in the ring of (Laurent) polynomials.

For a polynomial $p(t,u)=\sum p_{ij}t^iu^j$ its Newton polygon $N_p \subset \Z^2$ is defined as the convex hull of all $(i,j) \in \Z^2$ such that $p_{ij} \neq 0$. By the conditions of the lemma, $g(t,u)$ is not a monomial, thus $N_g$ has a non-trivial edge $v=(v_1, v_2)$. We will use this edge to find the number $k>1$ which divides all the necessary $(b-a)$ or $(c-b)$.

Suppose that $\det\nolimits_{a,b,c}+f_{a,b,c} \in S$. It is divisible by $GCD(S)$ and, moreover, by $g$. The polynomial $\det\nolimits_{a,b,c}+f_{a,b,c}$ is divisible by $\det\nolimits_{0,1,2}$, while $g$ is coprime with it, thus $(\det\nolimits_{a,b,c}+f_{a,b,c}) / \det\nolimits_{0,1,2}$ is divisible by $g$, that is
$$(\det\nolimits_{a,b,c}+f_{a,b,c})/\det\nolimits_{0,1,2} = g \cdot h,$$
where $h(t,u)$ is a (Laurent) polynomial. By the properties of the Newton polygon holds
$$N_{\det\nolimits_{a,b,c}+f_{a,b,c}} = N_g + N_h + N_{\det\nolimits_{0,1,2}},$$
where $+$ is the Minkowski sum $A+B=\{a+b: a\in A, b\in B\}$.

The Minkowski sum of two convex polygons in the plane can be computed as follows. Let us orientate the edges of each polygon counterclockwise. Then the edges of the Minkowski sum are sums of the parallel edges of the two convex polygons going in the same direction. Thus the edge $v$ is a part of some edge of $N_{\det\nolimits_{a,b,c}+f_{a,b,c}}$.

The polynomial
$$\det\nolimits_{a,b,c}(t,u)+f_{a,b,c}(t,u)=t^au^b-t^bu^a+t^bu^c-t^cu^b+t^cu^a-t^au^c+f_{a,b,c}(t,u)$$
has Newton polytope as shown on Figure \ref{fignewtdet}, and its edges oriented counterclockwise are proportional with positive coefficients to the vectors $(1,0)$, $(0,1)$, $(-1,1)$, $(-1,0)$, $(0,-1)$ and $(1,-1)$. 

\begin{figure}
    \includegraphics[width=8.0cm]{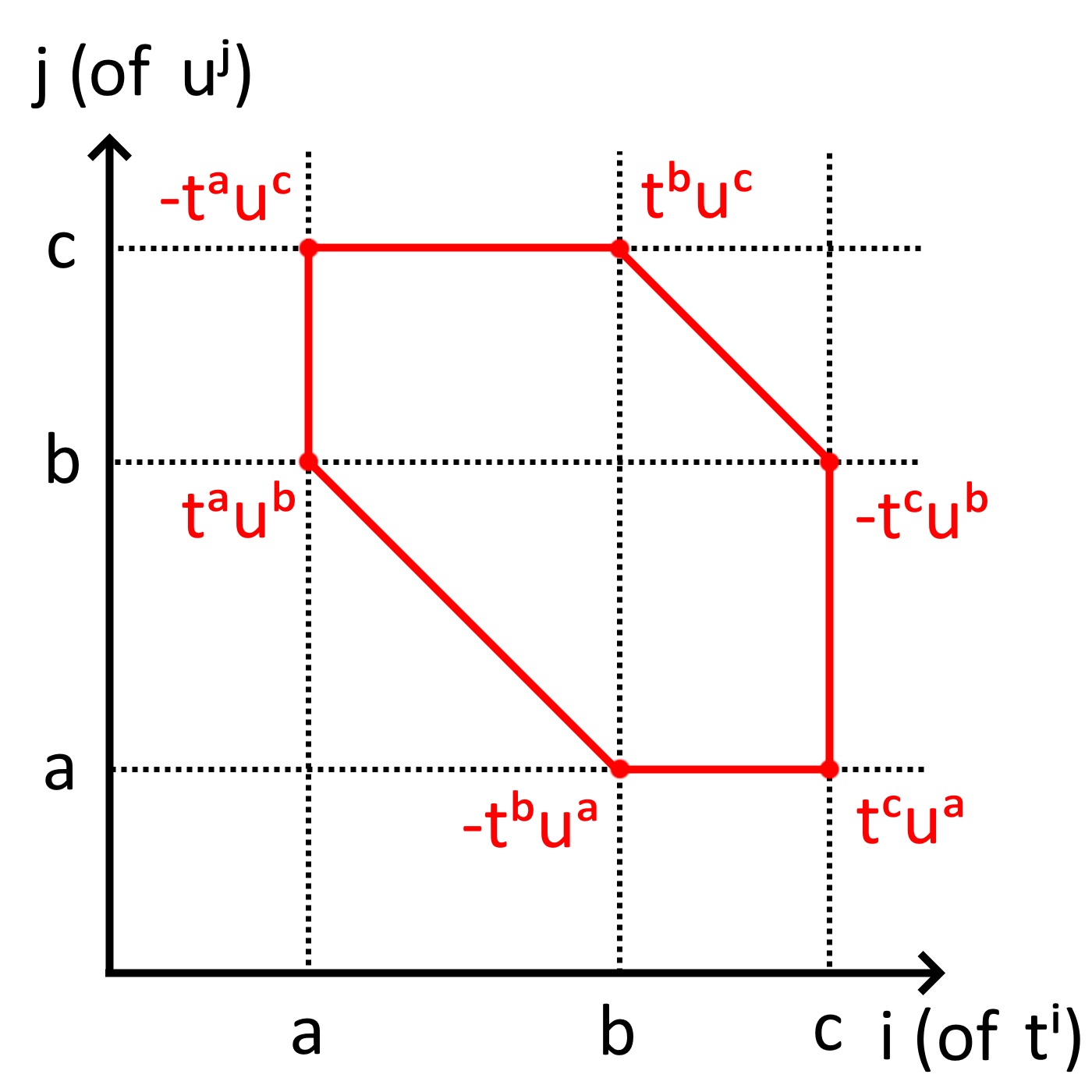}
    \caption{The Newton polytope of $\det\nolimits_{a,b,c}(t,u)$ with $a < b < c$}
    \label{fignewtdet}
\end{figure}

Let us now consider a weighted degree function $\deg_v(t^iu^j)=v_2 \cdot i - v_1 \cdot j$, for which the largest monomials of $N_g$ are the monomials that lie in the edge $v$. For a polynomial $p(t,u)=\sum p_{ij}t^iu^j$ let its degree with respect to $v$ be $\deg_v(p)=\max_{p_{ij} \neq 0} \deg_v(t^iu^j)$, and let its highest part with respect to $v$ be
$$p^v(t,u)=\sum_{\deg_v(t^iu^j)=\deg_v(p)} p_{ij}t^iu^j.$$
In particular, $g^v(t,u)=\sum_{(i,j) \in v} g_{ij}t^iu^j$, that is the highest part of $g(t,u)$ consists exactly of the monomials of $v$.

The highest part is multiplicative, thus
$$(\det\nolimits_{a,b,c}(t,u)+f_{a,b,c}(t,u))^v = g^v(t,u) \cdot h_1^v(t,u) \cdot \det\nolimits_{0,1,2}^v(t,u).$$

If $v$ is proportional with positive coefficient to $(1,0)$, then $$(\det\nolimits_{a,b,c}(t,u)+f_{a,b,c}(t,u))^v=-t^bu^a+t^cu^a=t^bu^a(t^{c-b}-1) \text{ and } \det\nolimits^v_{0,1,2}(t,u)=tu(t-1),$$
thus we have that $(t^{c-b}-1)/(t-1)$ is divisible by $g^v(t,u)$ (in the ring of Laurent polynomials).

We can now take the smallest natural number $k$ such that $(t^k-1)/(t-1)$ is divisible by $g^v(t,u)$ and denote it by $k$. The edge $v$ is non-empty, thus $g^v(t,u)$ is not a monomial and thus is non-invertible in the ring of Laurent polynomials, therefore $k \geqslant 2$.

The greatest common divisor of $(t^\alpha-1)$ and $(t^\beta-1)$ is $(t^{(\alpha, \beta)}-1)$. If for some $a$, $b$ and $c$ we have $\det\nolimits_{a,b,c}+f_{a,b,c} \in S$, then the polynomial $(t^{c-b}-1)/(t-1)$ is divisible by $g^v(t,u)$, thus $(t^{(c-b,k)}-1)/(t-1)$ is also divisible by $g^v(t,u)$, thus $(c-b,k) \geqslant k$, therefore $(c-b)$ is divisible by $k$. Consequently, $GCD(\{c-b  \,|\, \det\nolimits_{a,b,c}+f_{a,b,c} \in S\}) > 1$.

Now if $v$ is proportional with positive coefficient to $(-1,0)$, then $$(\det\nolimits_{a,b,c}(t,u)+f_{a,b,c}(t,u))^v=t^bu^c-t^au^c=t^au^c(t^{b-a}-1).$$
We can similarly define $k$ as the smallest natural number such that $(t^k-1)/(t-1)$ is divisible by $g^v(t,u)$, and similarly prove that if for some $a$, $b$ and $c$ we have $\det\nolimits_{a,b,c}+f_{a,b,c} \in S$, then $(b-a)$ is divisible by $k$.

The other four cases can be proven similarly. % TODO: more detail
\end{proof}

% TODO: something more clear in this section

\section{The strata S(1,2) and S(2,1) of the solution space for the filtration subset N(1,1,1)}\label{secstrats12forn111}

In this section we will study the stratum $S_{1,2}$ for $N(1,1,1)$, although without refering to it as such. This section is independent from the other parts of the text.

\begin{defin}
For a finite collection of integer numbers $n_1, n_2, \ldots, n_k$ let us define the function
$$\phi(n_1, n_2, \ldots, n_k) := (n_2-n_1, n_3-n_2, \ldots, n_k-n_{k-1})$$
to be the greatest common divisor of their differences.
\end{defin}

\begin{theor}\label{theorstratumranks21for111}
Fix $n \in \Z$. Let $x$ and $y$ be two $n$-th roots of unity, such that $x \neq 1$, $y \neq 1$, $x \neq y$. Let $B=\{ a, b, c, d, \ldots\} \subset \Z$ be a finite set of integer numbers. Consider a matrix 
$$M(B;x,y) = \begin{pmatrix}
1 & 1 & 1 & 1 & ...\\ 
x^a & x^b & x^c & x^d & ...\\
y^a & y^b & y^c & y^d & ...
\end{pmatrix},$$
where the columns are parametrized by $B$. If all the $3 \times 3$-minors of $M(B,x,y)$ are degenerate, then $B$ can be split into $B=B' \sqcup B''$ such that $\phi(B')$, $\phi(B'')$ and $n$ have a common divisor $k \geqslant 3$.
\end{theor}

This section is devoted to the proof of this theorem.

\subsection{Lemmas on 3x3-minors}

First, let us consider one $3 \times 3$-minor.

\begin{lemma}\label{lem3minor}
Let $x$ and $y$ be complex numbers of absolute value $1$, let $a$, $b$ and $c$ be integer numbers and consider a matrix
$$M=M(a,b,c;x,y) = \begin{pmatrix}
1 & 1 & 1\\ 
x^a & x^b & x^c \\
y^a & y^b & y^c
\end{pmatrix}$$
If this matrix is degenerate, then it has either two proportional rows or two proportional columns.
\end{lemma}
\begin{rem}
The determinant of this matrix is obviously related to a simplest Schur polynomial. In a forthcoming joint work with A. Voorhaar and A. Esterov, we shall use this relation in both directions: to apply known results on zero loci of Schur polynomials to a further study of singularities of resultants and discriminants, and to derive new results on zero loci of Schur polynomials from what we do here.
\end{rem}
\begin{proof}
Let us write the equation $\det M=0$ as
$$(x^a-x^b)(y^a-y^c)=(x^a-x^c)(y^a-y^b).$$
Suppose that $x^a=x^b$. Then either $x^a=x^c$ and $M$ has two proportional rows or $y^a=y^b$ and $M$ has two proportional (even equal) columns.

As the conditions of the lemma are symmetric in $a$, $b$ and $c$, and are symmetric in $x$ and $y$, we can now suppose that the numbers $x^a$, $x^b$, $x^c$ are pairwise different and that the numbers $y^a$, $y^b$, $y^c$ are also pairwise different. Thus we can write the condition above as
$$\frac{x^a-x^b}{x^a-x^c}=\frac{y^a-y^b}{y^a-y^c} \neq 0.$$

In particular, if we consider $x^a$, $x^b$ and $x^c$ as points in the plane, the following angles are equal: $\angle(x^b,x^a,x^c)=\angle(y^b,y^a,y^c)$. 

\begin{figure}
    \includegraphics[width=7.0cm]{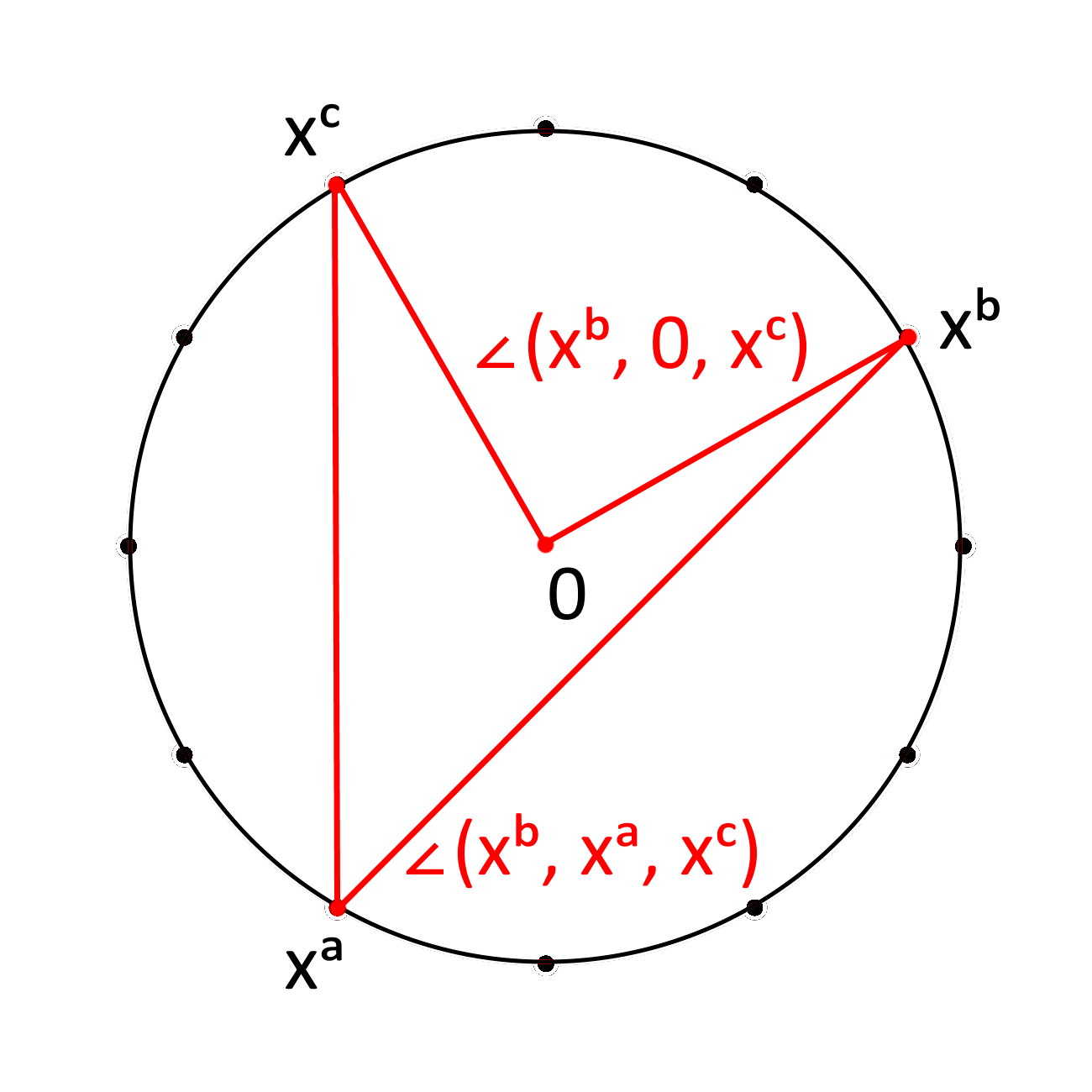}
    \caption{The points $x^a$, $x^b$ and $x^c$ on the unit circle}
    \label{figinscrangles}
\end{figure}

But the numbers $x$ and $y$ has absolute value $1$, thus the points $x^a$, $x^b$, $x^c$, $y^a$, $y^b$, $y^c$ lie on the unit circle and these angles are the inscribed angles on the arcs $(x^b, x^c)$ and $(y^b, y^c)$. An inscribed angle is half of the central angle on the same arc (see Figure \ref{figinscrangles}), thus
$$\angle(x^b,0,x^c)=2\angle(x^b,x^a,x^c)=2\angle(y^b,y^a,y^c)=\angle(y^b,0,y^c).$$

But $|x^b|=|x^c|=|y^b|=|y^c|=1$, thus $\angle(x^b,0,x^c)=\angle(y^b,0,y^c)$ implies $x^b/x^c=y^b/y^c$. As the conditions of the lemma are symmetric in $a$, $b$ and $c$, the matrix $M$ has two proportional rows.
\end{proof}

\iffalse
\begin{rem}\label{remlem3minor}
We can prove the second part of the lemma algerbraically. Indeed, let $\overline x$ be a complex conjugate to $x$. Then
$$\frac{x^a-x^b}{x^a-x^c} \Big/ \frac{\overline x^a-\overline x^b}{\overline x^a-\overline x^c}=\frac{x^a-x^b}{\overline x^a-\overline x^b} \cdot \frac{\overline x^a-\overline x^c}{x^a-x^c}=(-x^ax^b) \cdot (-x^{-a}x^{-c})=x^b/x^c,$$
thus $x^b/x^c=y^b/y^c$.
\end{rem}
\fi

Let us also prove a similar lemma for $J(2,1)$.

\begin{lemma}\label{lem3minor2}
Let $x$ be a complex number of absolute value $1$, let $a$, $b$ and $c$ be different integer numbers and consider a matrix
$$M=M(a,b,c;x,y) = \begin{pmatrix}
1 & 1 & 1\\ 
a & b & c \\
x^a & x^b & x^c
\end{pmatrix}$$
If this matrix is degenerate, then it has two equal rows, namely the rows 1 and 3.
\end{lemma}

\begin{proof}
Let us write the equation $\det M=0$ as
$$(a-b)(x^a-x^c)=(a-c)(x^a-x^b).$$
Suppose that $x^a=x^b$. Then $x^a=x^c$ and $M$ has two equal rows.

As the conditions of the lemma are symmetric in $a$, $b$ and $c$, we can now suppose that the numbers $x^a$, $x^b$, $x^c$ are pairwise different. Thus we can write the condition above as
$$\frac{x^a-x^b}{x^a-x^c}=\frac{a-b}{a-c} \neq 0.$$

The number on the right hand side is real, thus the plane vectors $x^a-x^b$ and $x^a-x^c$ are proportional. As the points $x^a$, $x^b$ and $x^c$ lie on the unit circle, it is possible only in the case $x^b=x^c$, which is a contradiction.
\end{proof}

\subsection{Proportionality of columns and rows}

\begin{lemma}\label{lemprop}
Consider a $3 \x n$-matrix $M$ with non-zero entries. Suppose that every its $3 \x 3$-minor has either two proportional rows or two proportional columns. Then either $M$ has two proportional rows or its columns can be divided into two groups of mutually proportional columns.
\end{lemma}

\begin{rem}
The lemma can be formulated in a more symmetric way. Let $M$ be a matrix with non-zero entries. Suppose that for every its $3 \x 3$-minor either its rows or its columns form at most two classes of equivalence with respect to proportionality. Then either the rows or the columns of the whole $M$ form at most two classes of equivalence with respect to proportionality. Moreover, one can prove it using the same logic even without the restriction that the dimensions of $M$ are $3\x n$.
\end{rem}

\begin{proof}
We will limit ourselves with the case when $M$ is $3 \x n$-matrix. Suppose that $M$ has three mutually non-proportional columns. Without loss of generality let them be the columns 1, 2 and 3. Consider the $3 \x 3$-minor given by these columns. By the conditions of the lemma it has two proportional rows. Without loss of generality let them be the rows 1 and 2. We would like to prove that the rows 1 and 2 of the whole $M$ are also proportional.

Consider the columns 1 and 2 which are triples of non-zero numbers. They are not proportional, while their 1-st and 2-nd entries are. Thus their 3-rd entries can not be proportional neither to their 1-st entries, nor to their 2-nd entries. Similarly, for the columns 1 and 3 (or 2 and 3) their 3-rd entries are not proportional neither to their 1-st entries, nor to their 2-nd entries.

Now consider the rows 1 and 2 of the whole $M$ and let us prove that they are proportional. Take their $k$-th elements. The column $k$ is proportional at most to one of the columns 1, 2 and 3. Without loss of generality suppose that it is not proportional neither to the column 1, nor to the column 2. Consider $3 \x 3$-minor given by the columns 1, 2 and $k$. Its columns are not proportional. The 1-st and 2-nd entries of its rows 1 and 3 are also not proportional, and the 1-st and 2-nd entries of its rows 2 and 3 are also not proportional. Thus by the conditions of the lemma the rows 1 and 2 of this minor should be proportional.

Consequently, the rows 1 and 2 of the whole $M$ are proportional.
\end{proof}

Now let us prove Theorem \ref{theorstratumranks21for111}.

\begin{proof}
By Lemma \ref{lem3minor}, the conditions of Lemma \ref{lemprop} are satisfied, thus either
$$M(B;x,y) = \begin{pmatrix}
1 & 1 & 1 & 1 & ...\\ 
x^a & x^b & x^c & x^d & ...\\
y^a & y^b & y^c & y^d & ...
\end{pmatrix}$$
has two proportional rows or its columns can be divided into two groups of mutually proportional columns.

It the rows 1 and 2 are proportional, then $x^a=x^b=x^c=x^d=\ldots$ and $y^a=y^b=y^c=y^d=\ldots$, thus (as in the proof of Lemma \ref{lemcork1simple} above) $x$ and $y$ are $\phi(B)$-th roots of unity. But they are also $n$-th roots of unity, thus they are $k=(\phi(B),n)$-th roots of unity. Moreover, as $x \neq y$ and they are different from 1, we have that $k \geqslant 3$. Thus we can define $B'=B$, $B''=\varnothing$. % TODO: check

If the columns, which are enumerated by $B$, can be divided into two groups of mutually proportional columns, then let the groups be $B'$ and $B''$. Then $x^{b_i}=x^{b_j}$ for $i, j \in B'$, thus $x$ is $\phi(B')$-th root of unity. Similarly, $x$ is $\phi(B'')$-th root of unity, and $y$ is $\phi(B')$-th and $\phi(B'')$-th root of unity. To sum it up, $x$ and $y$ are $k$-th roots of unity such that $\phi(B')$, $\phi(B'')$ and $n$ are divisible by $k$. Moreover, as $x \neq y$ and they are different from 1, we have that $k \geqslant 3$.
\end{proof}

\section{The projection of a spatial complete intersection curve}

In this section, we prove Theorem \ref{theorproj}.
%TODO: prove
%Suppose that $C=V(f_1, f_2)$ with $\supp f_i=A_i$. 
Let $f_i$ be a generic linear combination of monomials from a finite subset $A_i$ of the monomial lattice $\Z^3$.
Then we have a map
$$F:\CC^2 \to \C^{B_1} \times \C^{B_2}, \quad (x,y) \mapsto (f_1(x,y,\cdot), f_2(x,y,\cdot)),$$
and the sought closure $C$ of the image of $V(f_1,f_2)$ under the coordinate projection $\CC^3\to\CC^2$ equals the preimage of the resultant $R_B \subset \C^{B_1}\x\C^{B_2}$ under this map.

This observation allows to learn about singularities of $C$ from the singularities of $R_B$, using the following fact (see e.g. \cite{E08}, \cite{K16} or Theorem 10.7 in \cite{EL}).
\begin{utver}\label{ptransv}
For an algebraic set $V\subset\CC^n$ and generic linear combinations $h_i$ of monomials from non-empty finite sets $H_i\subset\Z^k,\,i=1,\ldots,k$, the intersection of $V$ and the zero locus $h_1=\cdots=h_k=0$ is a codimension $k$ subset in $V$, and this intersection is transversal at every smooth point of $V$. 

In particular, the intersection is empty if $\dim V<k$.
\end{utver}

Let us take any $x\in C=F^{-1}(R_B)$ and study $C$ in its neighborhood, depending on (1) what part of $R_B$ the pont $F(x)$ belongs to, and (2) what coordinate plane in $\C^{B_1} \times \C^{B_2}$ it belongs to.

More specifically, for every $I\subset B_1\sqcup B_2$, consider the coordinate plane $$\C^I:=\left\{\left(\sum_{b\in B_1\cap I}c_{1,b}t^b, \sum_{b\in B_2\cap I}c_{2,b}t^b\right)\right\}\subset\C^{B_1} \times \C^{B_2}.$$
Choose the unique $I$ such that $F(x)\in\CC^I$.

Recall that, by the part (i) of the main theorem \ref{theormain}, there exists a codimension 3 subset $\Sigma\subset \sing R_B$, such that at every point of $\sing R_B\setminus\Sigma$ the resultant is locally the union of two transversal smooth hypersurfaces, and, in particular, $\sing R_B$ is smooth at this point.

a) Assume $F(x)\in\CC^I\cap\Sigma$. This cannot happen by the dimension count, applying Proposition \ref{ptransv} to $\CC^n:=\CC^2\times\CC^I$ with coordinates $x_1,x_2,y_i,i\in I$, the set $V:=\CC^2\times(\CC^I\cap\Sigma)$, and the polynomials  $h_i(x,y)=f_i(x)-y_i,\,i\in I$ and $h_i=f_i,\,i\notin I$.

b) Otherwise, assume $F(x)\in\CC^I\cap\sing R_B$. Since $\sing(\CC^I\cap\sing R_B)$ has codimension at least 3, by the same dimension count, $F(x)$ is a regular point of $\CC^I\cap\sing R_B$. Moreover, at this point $F$ transversally intersects $\CC^I\cap\sing R_B$, applying Proposition \ref{ptransv} to the same $\CC^n$ and $h_i$'s and the set $V:=\CC^2\times(\CC^I\cap\sing R_B)$. Thus $F$ transversally intersects $\sing R_B$ at $F(x)\notin\Sigma$, thus $C=F^{-1}(R_B)$ is locally the union of two transversal smooth curves.

c) Otherwise, assume $F(x)\in\CC^I\cap R_B$. Then, by the same dimension count, $F(x)$ is a regular point of $\CC^I\cap R_B$. Moreover, at this point $F$ transversally intersects $\CC^I\cap R_B$, applying Proposition \ref{ptransv} to the same $\CC^n$ and $h_i$'s and the set $V:=\CC^2\times(\CC^I\cap R_B)$. Thus $F$ transversally intersects $ R_B$ at $F(x)\notin\sing R_B$, thus $C=F^{-1}(R_B)$ is locally a smooth curve.

\end{document}